\documentclass{ieeeaccessarxiv}

\usepackage{cite}
\usepackage{amsmath,amssymb,amsfonts,amsthm,mathtools}
\usepackage{graphicx}
\usepackage{booktabs}
\usepackage{array}
\usepackage{enumitem}
\usepackage{xurl}
\usepackage[hidelinks]{hyperref}
\hypersetup{
  pdftitle={Information-Theoretic Upper Bounds for Deterministic Noise in Zeroth-Order Convex Optimization},
  pdfauthor={Dmitry Pasechnyuk-Vilensky; Igor Pavlov; Martin Tak\'a\v{c}; Alexander Gasnikov},
  pdfsubject={Finite-budget noise tolerance in zeroth-order convex optimization},
  pdfkeywords={Derivative-free optimization, deterministic noise, maximum admissible noise, Moreau envelope, noisy oracle, zeroth-order convex optimization}
}

\allowdisplaybreaks

\theoremstyle{plain}
\newtheorem{theorem}{Theorem}[section]
\newtheorem{lemma}[theorem]{Lemma}
\newtheorem{proposition}[theorem]{Proposition}
\newtheorem{corollary}[theorem]{Corollary}
\theoremstyle{definition}
\newtheorem{definition}[theorem]{Definition}
\newtheorem{remark}[theorem]{Remark}

\newcommand{\R}{\mathbb{R}}
\newcommand{\Prob}{\mathbb{P}}
\newcommand{\E}{\mathbb{E}}
\newcommand{\eps}{\varepsilon}
\newcommand{\cF}{\mathcal{F}}
\newcommand{\prox}{\operatorname{prox}}
\newcommand{\Proj}{\operatorname{Proj}}
\newcommand{\spn}{\operatorname{span}}
\newcommand{\norm}[1]{\left\|#1\right\|_2}
\newcommand{\ip}[2]{\langle #1,#2\rangle}
\newcommand{\oo}{\mathsf O}
\newcommand{\Sph}{\mathbb S^{n-1}}

\makeatletter
\renewcommand\appendices{\par
  \def\theHsection{Appendix.\Alph{section}}%
  \setcounter{section}{-1}\setcounter{subsection}{0}%
  \setcounter{subsubsection}{0}\setcounter{paragraph}{0}%
  \def\thesection{\Alph{section}}%
  \def\thesectiondis{APPENDIX \Alph{section}.}%
  \refstepcounter{section}\setcounter{section}{0}}
\makeatother

\begin{document}
\history{}
\doi{}

\title{Information-Theoretic Upper Bounds for Deterministic Noise in Zeroth-Order Convex Optimization}

\author{\uppercase{Dmitry Pasechnyuk-Vilensky}\authorrefmark{1,2}, \uppercase{Igor Pavlov}\authorrefmark{3},
\uppercase{Martin Tak\'a\v{c}}\authorrefmark{1}, and
\uppercase{Alexander Gasnikov}\authorrefmark{2,4}}
\address[1]{Mohamed bin Zayed University of Artificial Intelligence, Abu Dhabi, UAE  (e-mail: dmivilensky1@gmail.com)}
\address[2]{Trusted AI Research Center, RAS, Moscow, Russia}
\address[3]{Moscow Independent Research Institute of Artificial Intelligence, Moscow, Russia}
\address[4]{Innopolis University, Innopolis, Russia}
\tfootnote{The study was supported by the Ministry of Economic Development of the Russian Federation (agreement with MIPT No. 139-15-2025-013, dated June 20, 2025, IGK 000000C313925P4B0002).}
\markboth
{Pasechnyuk-Vilensky \headeretal: Information-Theoretic Upper Bounds for Deterministic Noise in Zeroth-Order Convex Optimization}
{Pasechnyuk-Vilensky \headeretal: Information-Theoretic Upper Bounds for Deterministic Noise in Zeroth-Order Convex Optimization}
\corresp{Corresponding author: Dmitry Pasechnyuk-Vilensky.}
\begin{abstract}
We study zeroth-order convex optimization on Euclidean balls with function values corrupted by a fixed, uniformly bounded deterministic perturbation. For a query budget $T$, the maximum admissible level of noise (MALN) is the supremum of noise levels for which an $\eps$-accurate solution can be found with probability $1-\beta$. We construct an explicit hard family: the support function of a spherical belt combined with an affine branch along a hidden random direction. It yields finite-budget upper bounds on the MALN for Lipschitz, strongly convex, uniformly convex, smooth, H\"older, and smooth strongly convex classes. For Lipschitz objectives the bound has order $\eps^2\sqrt\ell/(\sqrt nMR)+\eps\ell/n$, where $\ell=\log(2(T+1)/(1-\beta))$ and $n-1\ge32\ell$; this is the scale of Li and Risteski with explicit logarithmic factors. A tangential-gradient estimate for the Moreau envelope yields the smooth scale $\eps^{3/2}/(\sqrt n\sqrt L\,R)$ and, for $L-\mu\ge5\mu$, the smooth strongly convex scale $\eps\sqrt{\mu/((L-\mu)n)}$. These bounds are tight: on explicit parameter ranges and subclasses, polynomial-query two-point methods tolerate noise of the order of the terms of order $n^{-1/2}$, so that the MALN is determined up to $O(\sqrt\ell)$ and absolute constants whenever these terms dominate. At $L=\mu$, simplex interpolation with $n+1$ queries tolerates noise $R\sqrt{2\mu\eps/n}$, which is optimal up to $O(\sqrt\ell)$ for $\eps\le\mu R^2/3$, while with at most $n$ exact queries no randomized algorithm guarantees success probability greater than $1/2$ uniformly over the class. We also quantify noise tolerance near this endpoint and illustrate the hiding mechanism and the breakdown of a two-point method numerically.
\end{abstract}

\begin{keywords}
Derivative-free optimization, deterministic noise, maximum admissible noise, Moreau envelope, noisy oracle, zeroth-order convex optimization.
\end{keywords}

\titlepgskip=-15pt
\maketitle

\section{Introduction}
\PARstart{Z}{eroth-order} optimization concerns settings in which an algorithm can query objective values but has no direct access to gradients or subgradients. This model is central in derivative-free optimization, simulation-based design, black-box machine learning, and bandit convex optimization; see, for example, \cite{NemirovskiYudin1983,ConnScheinbergVicente2009,Spall2003,LarsonMenickellyWild2019,FlaxmanKalaiMcMahan2005,BubeckLeeEldan2017}. In many applications the returned function values are not exact. The error may come from finite precision, numerical inner solvers, simulation error, adversarial perturbations, robust worst-case modeling, or inexact oracle access \cite{DevolderGlineurNesterov2014,dAspremont2008,BerahasCaoChoromanskiScheinberg2022}. Systematic errors of this kind are deterministic functions of the query point, and repeated evaluation returns the same value, so averaging does not remove them.

We consider a deterministic noisy oracle. For a function $f$ on the ball $B_R\subset\R^n$, a query at $x$ returns a fixed value within distance $\delta$ of $f(x)$. The error is uniformly bounded but otherwise arbitrary. Besides the usual oracle-complexity question, this raises a noise-tolerance question: how large may $\delta$ be before a prescribed number of queries no longer guarantees an $\eps$-minimizer? We call this critical scale the maximum admissible level of noise, or MALN. An upper bound on the MALN is an impossibility statement; a lower bound on the MALN is an algorithmic guarantee.

We restrict attention to convex objectives. For Lipschitz convex functions on a ball, exact values permit optimization with a number of queries polynomial in the dimension and the inverse normalized accuracy \cite{NemirovskiYudin1983,Protasov1996}, whereas the worst-case exact-value complexity for general Lipschitz nonconvex objectives is exponential in the dimension \cite{NemirovskiYudin1983}. Noise introduces an additional obstruction, whose effect must be assessed at a specified query budget. An impossibility result also applies to any larger class containing the hard instances used in its proof. In particular, convex instances yield bounds for quasiconvex and weakly convex classes; extension to a Polyak--{\L}ojasiewicz class requires verifying that the chosen instances satisfy its defining condition.

\subsection{Contributions}
The results are organized in five blocks, summarized in Table~\ref{tab:summary} and Fig.~\ref{fig:roadmap}.
\begin{enumerate}[leftmargin=*]
\item {A finite-budget hiding principle and an explicit hard family.} Lemma~\ref{lem:hiding} reduces the upper bounds on the MALN, except the query threshold of Theorem~\ref{thm:main-endpoint}(iii), to four verifiable properties of a random family and a fixed common oracle. The family is built from the support function of a spherical belt and a hidden affine branch (Lemma~\ref{lem:belt}). All statements are finite-dimensional, with explicit constants and an explicit condition relating $n$, the budget $T$, and the failure probability $\beta$.
\item {Lipschitz and strongly convex bounds.} Theorems~\ref{thm:main-lip} and \ref{thm:main-lipsc} give the Lipschitz scale of Li and Risteski and its strongly convex analogue with the factors $\sqrt\ell$ and $\ell$ in the two terms, where $\ell=\log(2(T+1)/(1-\beta))$. Corollary~\ref{cor:uc} extends the strongly convex bound to uniformly convex objectives, and Proposition~\ref{prop:lipsc-end} shows that the Lipschitz strongly convex class becomes trivial when $M-\mu R\le\sqrt{2\mu\eps}$.
\item {Smooth bounds through a family-specific Moreau lemma.} Lemma~\ref{lem:moreau} bounds the discrepancy between the Moreau envelope of the hard family and a common radial function. Theorems~\ref{thm:main-smooth} and \ref{thm:main-smsc} and Corollary~\ref{cor:holder} give the smooth convex, smooth strongly convex, and H\"older scales. Proposition~\ref{prop:barrier} shows that a transfer from the Lipschitz class through a universal convolution cannot reach the smooth scale.
\item {Endpoint $L=\mu$.} Theorem~\ref{thm:main-endpoint} determines the endpoint MALN up to $O(\sqrt\ell)$ for $T\ge n+1$, $n-1\ge32\ell$, and $\eps\le\mu R^2/3$, identifies the exact query threshold $n+1$ for randomized algorithms, and Proposition~\ref{prop:near-end} extends the positive result to a quantitative neighborhood of the endpoint.
\item {Tightness and budgets.} Propositions~\ref{prop:ach-lip}--\ref{prop:ach-smsc} give, in the present model, polynomial-query two-point methods, the first two with the orders of \cite{DvinskikhTominin2022,AlashqarGasnikovDvinskikhLobanov2023}, whose noise tolerance has the order of the terms of order $n^{-1/2}$ of the upper bounds. By Corollary~\ref{cor:two-sided}, with $\theta=\theta_*$ and bounded $\kappa$, the upper bounds match the lower bounds up to $O(\sqrt\ell)$ and absolute constants whenever the terms linear in $\theta$ dominate, in the stated parameter and budget ranges: for Theorems~\ref{thm:main-lip} and \ref{thm:main-lipsc} (the latter for $M\ge2\mu R$ and $M-\mu R\ge4\sqrt{\mu\eps}$), for Theorem~\ref{thm:main-smooth} and Corollary~\ref{cor:holder} on the corresponding Lipschitz subclasses, and for Theorem~\ref{thm:main-smsc}(b) on the subclass of functions whose minimizer lies at distance at least $\sqrt{\eps/L}$ from the boundary. The terms $\eps\theta_*^2=\Theta(\eps\ell/n)$ correspond to the scale $\eps/n$ of the simulated-annealing method of \cite{BelloniLiangNarayananRakhlin2015}. Corollary~\ref{cor:budget} converts the Lipschitz upper bound into an exponential lower bound on the number of queries at a given noise level.
\end{enumerate}

\subsection{Related Work}
Singer and Vondr\'ak \cite{SingerVondrak2015} established information-theoretic lower bounds for convex optimization with erroneous value oracles. Belloni, Liang, Narayanan, and Rakhlin \cite{BelloniLiangNarayananRakhlin2015} gave polynomial-query methods for functions uniformly $O(\eps/n)$-close to convex ones through simulated annealing, which is the origin of the universal branch $\eps/n$. Li and Risteski \cite{LiRisteski2016} identified, up to logarithmic factors, the threshold $\max\{\eps^2/(\sqrt n MR),\eps/n\}$ for $M$-Lipschitz convex functions: below the threshold divided by an absolute constant a polynomial-time method succeeds on well-conditioned bodies, including the ball \cite[Thm.~3.2]{LiRisteski2016}, and for every $c\ge1$, every $n\ge d_c$, and every algorithm there are an accuracy $\eps<1/64$ and such a function on the ball of diameter one, with noise level at least the threshold multiplied by $(13c\log(n/\eps))^2$, on which the algorithm, run for time $o((n/\eps)^c)$, succeeds with probability at most $1/2$ \cite[Thm.~3.1 and Sec.~4]{LiRisteski2016}. Their approximately convex model uses the same perturbations as the present deterministic model: a Borel function $F$ with $\sup|F-f|\le\delta$ for a convex $f$ is a fixed oracle in $\oo_\delta(f)$. Their accuracy criterion $F(\tilde x)\le\min F+\eps$ refers to $F$; since $|F-f|\le\delta$, it differs from the criterion $f(\tilde x)\le\min f+\eps$ used here by at most $2\delta$ in the accuracy. In the normalization $M=1$, $R=1/2$ and for budgets $T=(n/\eps)^c$, Theorem~\ref{thm:main-lip} carries the factor $\sqrt\ell=O(\sqrt{c\log(n/\eps)})$ in the term of order $n^{-1/2}$ and $\ell=O(c\log(n/\eps))$ in the term $\eps/n$, in place of $(13c\log(n/\eps))^2$, and applies to every fixed budget $T$. The hard family of Section~\ref{sec:main-results} follows their principle of a hidden linear branch in a random double cone, with a different radial part that makes the common oracle exact on visible points beyond an explicit radius (Remark~\ref{rem:li-risteski}).

Algorithmic noise tolerance of randomized gradient-free methods is analyzed in \cite{NesterovSpokoiny2017,Gasnikov2024,GorbunovDvurechenskyGasnikov2022,DvurechenskyGorbunovGasnikov2021,BogolubskyEtAl2016,Kornilov2023GradientFree,Kornilov2023Heavy}, and two-point estimators with optimal dependence on the dimension are studied in \cite{DuchiJordanWainwrightWibisono2015,Shamir2017}. For bounded deterministic noise, the admissible level of two-point methods is studied in \cite{BayandinaGasnikovLagunovskaya2018,DvinskikhTominin2022,AlashqarGasnikovDvinskikhLobanov2023}; see also the survey \cite{Gasnikov2024}. For Lipschitz convex objectives on a set of diameter $D$, Dvinskikh, Tominin, Tominin, and Gasnikov \cite{DvinskikhTominin2022} obtained the level $\eps^2/(\sqrt nMD)$ with $O(nM^2D^2/\eps^2)$ queries and observed, from the bound of \cite{LiRisteski2016}, that for $\eps\gtrsim MD/\sqrt n$ this level is the largest admissible one up to a logarithmic factor; the same level is derived for $l_1$- and $l_2$-randomized smoothing in \cite{AlashqarGasnikovDvinskikhLobanov2023}, and the level $\mu_r^{1/r}\eps^{2-1/r}/(M\sqrt n)$ is indicated under an $r$-growth condition in \cite{DvinskikhTominin2022}. For objectives of higher-order smoothness, kernel-based methods tolerate noise of order $\eps^{3/2}/\sqrt n$ in the dependence on $\eps$ and $n$ under bounded-variance, possibly adversarial, noise \cite{LobanovBashirovGasnikov2024}, and the maximal noise level is analyzed under the Polyak--{\L}ojasiewicz condition in \cite{GasnikovLobanovStonyakin2024}. Under the direction-independence and bounded-second-moment assumptions of \cite{Lobanov2023Stochastic}, the stochastic perturbation contributes no bias to the two-point estimator; a level of order $\eps/\sqrt n$ preserves the stated oracle-complexity order. This is not a universal impossibility threshold for stochastic noise. Lobanov and Gasnikov \cite{LobanovGasnikov2024LipschitzNoise} instead bound the Lipschitz constant of the deterministic perturbation, not its uniform amplitude. Their noise parameter therefore differs from $\delta$ in the present model. Higher-order smoothness in zeroth-order optimization, including adversarial noise, is exploited in \cite{BachPerchet2016,AkhavanPontilTsybakov2020,AkhavanPontilTsybakov2021,AkhavanChzhenPontilTsybakov2024,BychkovEtAl2024}. These works concern a variety of noise models; Table~\ref{tab:models} records the distinctions relevant here. For exact function values, Protasov \cite{Protasov1996} gave a deterministic method whose number of queries is quadratic in the dimension up to logarithmic factors. For stochastic evaluations, Jamieson, Nowak, and Recht \cite{JamiesonNowakRecht2012} and Shamir \cite{Shamir2013} proved lower bounds in derivative-free and bandit settings. For exact values on the Euclidean ball, Kerger \cite{Kerger2026} proved that deterministic algorithms need $\Omega(n^2/\log n)$ queries at an accuracy of order $MR/\sqrt n$, and Zhang, Zhang, Qi, and Lin \cite{ZhangZhangQiLin2026} proved the lower bound $\Omega(n\min\{n,M^2R^2/\eps^2\}/\log\min\{n,M^2R^2/\eps^2\})$ for randomized algorithms; the hard instances are maxima of linear functions and random support functions, respectively. Inexact first-order oracles are treated in \cite{DevolderGlineurNesterov2014,dAspremont2008}, universal methods for H\"older continuous gradients in \cite{Nesterov2015Universal}, restart schemes for uniformly convex objectives in \cite{IouditskiNesterov2014}, and randomized smoothing in \cite{DuchiBartlettWainwright2012}.

\section{Model, Notation, and Function Classes}\label{sec:model}

For $R>0$ let $B_R=B_2^n(R)=\{x\in\R^n:\norm{x}\le R\}$, and let $\Sph$ be the unit sphere. Throughout, $n\ge2$, $T$ is a nonnegative integer, $0<\beta<1$, the quantities $\eps,R,M,L,\mu$ are positive, and $Q(x):=\frac\mu2\norm{x}^2$. For a nonempty closed convex set $C\subset\R^n$, $\Proj_C$ denotes the Euclidean projection onto $C$. The norm is Euclidean for vectors and the induced operator norm for matrices. The class parameters, the target accuracy, the failure probability, the noise bound, and the budget are known to the algorithm; the objective and its oracle are not.

\begin{definition}[Deterministic noisy zeroth-order oracle]
Let $f:B_R\to\R$ and $\delta\ge0$. Define
\begin{equation*}
\oo_\delta(f)
:=
\{O:B_R\to\R\ \text{Borel}:\ \sup_{x\in B_R}|O(x)-f(x)|\le\delta\}.
\end{equation*}
\end{definition}

A {randomized adaptive algorithm with $T$ queries} consists of a standard Borel probability space $(\Omega,\Prob_A)$ and Borel maps $\chi_t:\Omega\times\R^{t-1}\to B_R$, $t=1,\ldots,T+1$. Against an oracle $O$ it produces the queries $x_t=\chi_t(\omega,y_1,\ldots,y_{t-1})$ with answers $y_t=O(x_t)$, $t\le T$, and the output $A^O(\omega)=\chi_{T+1}(\omega,y_1,\ldots,y_T)$. An algorithm with fewer queries is the special case in which the maps ignore some answers. Deterministic algorithms correspond to a one-point $\Omega$.

\begin{definition}[Solvability and MALN]\label{def:maln}
An algorithm $A$ with $T$ queries {$(\eps,\beta,\delta,T)$-solves} a class $\cF$ of functions on $B_R$ if, for every $f\in\cF$ and every $O\in\oo_\delta(f)$,
\begin{equation*}
\Prob_A\Bigl\{
f(A^O)-\inf_{B_R} f\le\eps
\Bigr\}
\ge 1-\beta.
\end{equation*}
The class $\cF$ is {$(\eps,\beta,\delta,T)$-solvable} if some algorithm with $T$ queries $(\eps,\beta,\delta,T)$-solves it. The maximum admissible level of noise is
\begin{equation*}
\Delta_{\cF}(\eps,\beta,T)
:=
\sup\{\delta\ge0:\ \cF\text{ is }(\eps,\beta,\delta,T)\text{-solvable}\},
\end{equation*}
with $\sup\varnothing:=0$. Its values lie in $[0,+\infty]$.
\end{definition}

Since $\oo_{\delta'}(f)\subset\oo_\delta(f)$ for $\delta'\le\delta$, an algorithm that $(\eps,\beta,\delta,T)$-solves $\cF$ also $(\eps,\beta,\delta',T)$-solves it. Consequently, if no algorithm $(\eps,\beta,D,T)$-solves $\cF$, then $\Delta_\cF(\eps,\beta,T)\le D$. By convention, $\Delta_\cF=0$ does not by itself imply solvability with exact values: it also covers the case in which the set in Definition~\ref{def:maln} is empty. A supremum need not be attained.

We use the following classes of functions $f:B_R\to\R$:
\begin{itemize}[leftmargin=*]
\item $\cF_{\mathrm{Lip}}^M(B_R)$: $f$ convex and $M$-Lipschitz;
\item $\cF_{\mathrm{Sm}}^L(B_R)$: $f$ convex, the restriction to $B_R$ of a function differentiable on an open neighborhood of $B_R$, and $\nabla f$ is $L$-Lipschitz on $B_R$;
\item $\cF_{\mathrm{H}}^{\nu,L_\nu}(B_R)$, $\nu\in(0,1]$: $f$ convex, the restriction to $B_R$ of a function differentiable on an open neighborhood of $B_R$, and $\norm{\nabla f(x)-\nabla f(y)}\le L_\nu\norm{x-y}^\nu$ on $B_R$;
\item $\cF_{\mathrm{LipSC}}^{M,\mu}(B_R)$ and $\cF_{\mathrm{SmSC}}^{L,\mu}(B_R)$: the members of $\cF_{\mathrm{Lip}}^M(B_R)$ and $\cF_{\mathrm{Sm}}^L(B_R)$, respectively, for which $f-Q$ is convex on $B_R$;
\item $\cF_{\mathrm{UC}}^{M,\mu_p,p}(B_R)$, $p\ge2$: the members of $\cF_{\mathrm{Lip}}^M(B_R)$ such that
$f(y)\ge f(x)+\ip{s}{y-x}+\frac{\mu_p}{p}\norm{y-x}^p$ for all $x,y\in B_R$ and all $s\in\partial f(x)$, where $\partial f(x)=\{s:f(y)\ge f(x)+\ip s{y-x}\ \forall y\in B_R\}$.
\end{itemize}
Intersections such as $\cF_{\mathrm{Lip}}^M\cap\cF_{\mathrm{Sm}}^L$ are taken on the same ball; $\cF_{\mathrm{H}}^{1,L}=\cF_{\mathrm{Sm}}^L$. For the differentiable classes, the gradient at a boundary point of $B_R$ is determined by $f$ itself, because the directional derivatives along the directions pointing into $B_R$ span $\R^n$. The subdifferential in the uniform-convexity definition is relative to $B_R$ and may contain unbounded outward normal vectors at boundary points.

\begin{definition}[Hiding parameter]\label{def:hiding}
Let $T\ge0$ be an integer and $0<\beta<1/2$. A number $\theta\in(0,1/4]$ is {$(n,T,\beta)$-hiding} if
\begin{equation}\label{eq:hiding}
2(T+1)(1-\theta^2)^{(n-1)/2}<1-\beta .
\end{equation}
\end{definition}
The left side of \eqref{eq:hiding} decreases in $\theta$, and $(1-\theta^2)^{(n-1)/2}<e^{-(n-1)\theta^2/2}$. Hence every $\theta\in[\theta_*,1/4]$ is hiding, where
\begin{equation}\label{eq:theta-star}
\theta_*=\theta_*(n,T,\beta):=\sqrt{\frac{2\ell}{n-1}},\qquad \ell:=\log\frac{2(T+1)}{1-\beta},
\end{equation}
provided $\theta_*\le1/4$, that is, $\ell\le(n-1)/32$. For budgets polynomial in $n$ this gives $\theta_*=O(\sqrt{\log n/n})$. For a hidden direction $w\in\Sph$ we write
\[
A_w:=\{x\in\R^n\setminus\{0\}:\ |\ip wx|>\theta\norm x\},
\]
the open double cone of half-angle $\arccos\theta$ around $\pm w$; points outside $A_w$ are called {visible}. The hiding parameter $\theta$ is fixed in each statement and suppressed from the notation $A_w$.

\begin{lemma}[Finite-budget hiding]\label{lem:hiding}
Let $\cF$ be a class of Borel functions on $B_R$, let $(W,\Prob_W)$ be a standard Borel probability space, and let $D\ge0$ and $p\in[0,1]$. Suppose that for every $w\in W$ there are $f_w\in\cF$, a Borel oracle $O_w$, and a Borel set $E_w\subset B_R$, and that there is one Borel function $H:B_R\to\R$ independent of $w$, such that $(w,x)\mapsto\mathbf 1\{x\in E_w\}$ is jointly Borel and
\begin{enumerate}[label=(\alph*),leftmargin=*]
\item $O_w\in\oo_D(f_w)$;
\item $O_w=H$ on $B_R\setminus E_w$;
\item every $x\in B_R$ with $f_w(x)-\inf_{B_R}f_w\le\eps$ belongs to $E_w$;
\item $\Prob_W\{x\in E_w\}\le p$ for every $x\in B_R$.
\end{enumerate}
If $(T+1)p<1-\beta$, then no randomized algorithm with $T$ queries $(\eps,\beta,D,T)$-solves $\cF$. In particular $\Delta_\cF(\eps,\beta,T)\le D$. For $W=\Sph$ with the uniform law and $E_w=A_w\cap B_R$, (d) holds with $p=\min\{1,2(1-\theta^2)^{(n-1)/2}\}$, and the condition $(T+1)p<1-\beta$ is \eqref{eq:hiding}.
\end{lemma}

The proof in Appendix~\ref{app:lip} fixes the internal randomness, runs the algorithm against the common function $H$, and applies a union bound along the resulting {reference transcript}.

\subsection{The Hard Family}
For $w\in\Sph$ and $\theta\in(0,1/4]$ let
\begin{gather*}
K_w:=\{v\in\R^n:\norm v\le1,\ |\ip wv|\le\theta\},\\
\sigma_w(x):=\max_{v\in K_w}\ip vx
\end{gather*}
be the spherical belt and its support function. For parameters $0<a\le m$ and $\rho>0$ define, on all of $\R^n$,
\begin{equation}\label{eq:belt}
g_w(x):=\max\{m\sigma_w(x)-b,\ -a\ip wx\},\quad b:=\rho(m\theta+a),
\end{equation}
together with the radial comparison function $G(x):=\max\{m\norm x-b,0\}$ and
\begin{equation}\label{eq:D0}
\begin{gathered}
r_0:=\frac{b}{m-a\theta},\qquad D_0:=a\theta r_0,\\
D_1(\lambda):=a(r_0+\lambda m)\arcsin\theta .
\end{gathered}
\end{equation}
Lemma~\ref{lem:belt} in Appendix~\ref{app:lip} shows that $g_w$ is convex and $m$-Lipschitz, that its global minimum $-a\rho$ is attained at $\rho w$, and that $|g_w-G|\le D_0$ on visible points, with equality $g_w=G$ at visible points of norm at least $r_0$. With $q:=a/m$,
\begin{equation}\label{eq:r0-ratio}
\frac{r_0}{\rho}=\frac{\theta+q}{1-q\theta}.
\end{equation}

\section{Main Results}\label{sec:main-results}

\subsection{Outline of Results}
\begin{figure}[t]
\centering
\includegraphics[width=\columnwidth]{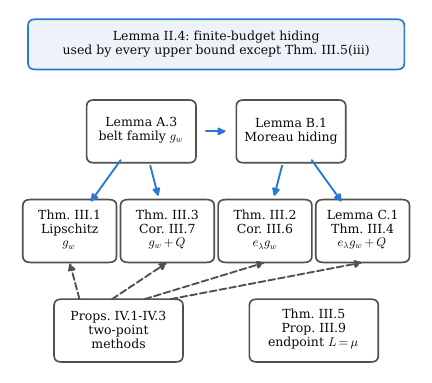}
\caption{Dependencies between the results. Lemma~\ref{lem:hiding} enters every upper bound except Theorem~\ref{thm:main-endpoint}(iii); solid arrows show the family-specific dependencies, and dashed arrows the comparisons of Corollary~\ref{cor:two-sided}.}
\label{fig:roadmap}
\end{figure}

Except for the exact-value query threshold in Theorem~\ref{thm:main-endpoint}(iii), all upper bounds on the MALN follow one scheme. A hard function $f_w$ is chosen with a hidden direction $w$ uniform on $\Sph$; every $\eps$-minimizer of $f_w$ lies in the double cone $A_w$; a fixed oracle agrees with a common function outside $A_w$; and Lemma~\ref{lem:hiding} converts the uniform discrepancy between $f_w$ and the common function into an upper bound on the MALN. The four main classes differ only in the choice of $f_w$: the belt $g_w$ (Lipschitz), $g_w+Q$ (strongly convex), the Moreau envelope $e_\lambda g_w(x)=\min_y\{g_w(y)+\norm{x-y}^2/(2\lambda)\}$ (smooth), and $e_\lambda g_w+Q$ (smooth strongly convex). Fig.~\ref{fig:roadmap} records the dependencies. The endpoint result, Theorem~\ref{thm:main-endpoint}, combines quadratic rigidity, simplex interpolation, a Gaussian reflection argument, and the same hiding lemma.

\subsection{Bounds on the MALN}
For Theorems~\ref{thm:main-lip}--\ref{thm:main-smsc} and Corollaries~\ref{cor:holder} and \ref{cor:uc}, assume $0<\beta<1/2$, $T\ge0$, and that $\theta$ is $(n,T,\beta)$-hiding. Theorem~\ref{thm:main-endpoint} states its own assumptions; its parts (i)--(iii) require no hiding condition.

\begin{theorem}[Lipschitz convex]\label{thm:main-lip}
If $0<\eps\le MR/8$, then no algorithm $(\eps,\beta,D,T)$-solves $\cF_{\mathrm{Lip}}^M(B_R)$ with $D=\frac43\eps\theta^2+\frac53\frac{\eps^2\theta}{MR}$. Consequently,
\begin{equation*}
\Delta_{\cF_{\mathrm{Lip}}^M(B_R)}(\eps,\beta,T)
\le
\frac43\eps\theta^2+\frac53\,\frac{\eps^2\theta}{MR}.
\end{equation*}
If $\ell\le(n-1)/32$, then $\theta=\theta_*$ from \eqref{eq:theta-star} is admissible and the bound is $\frac83\frac{\eps\ell}{n-1}+\frac{5\sqrt2}{3}\frac{\eps^2\sqrt\ell}{\sqrt{n-1}\,MR}$.
\end{theorem}

\begin{theorem}[Smooth convex]\label{thm:main-smooth}
If $0<\eps\le LR^2/16$, then
\begin{equation*}
\Delta_{\cF_{\mathrm{Sm}}^L(B_R)}(\eps,\beta,T)
\le
\frac{12}5\eps\theta^2+7\,\frac{\eps^{3/2}\theta}{\sqrt L\,R},
\end{equation*}
and the same bound holds for $\cF_{\mathrm{Sm}}^L(B_R)\cap\cF_{\mathrm{Lip}}^M(B_R)$ whenever $M\ge\sqrt{L\eps}$.
\end{theorem}

\begin{theorem}[Lipschitz strongly convex]\label{thm:main-lipsc}
If $0<\eps\le\mu R^2$ and $M-\mu R\ge4\sqrt{\mu\eps}$, then
\begin{equation*}
\Delta_{\cF_{\mathrm{LipSC}}^{M,\mu}(B_R)}(\eps,\beta,T)
\le
\frac{16}{7}\Bigl(\eps\theta^2+\frac{2\sqrt\mu\,\eps^{3/2}\theta}{M-\mu R}\Bigr).
\end{equation*}
\end{theorem}

\begin{theorem}[Smooth strongly convex]\label{thm:main-smsc}
Let $L_0:=L-\mu>0$.
\begin{enumerate}[label=(\alph*),leftmargin=*]
\item If $L_0\ge36\mu$ and $0<\eps\le\mu R^2$, then
\begin{equation*}
\Delta_{\cF_{\mathrm{SmSC}}^{L,\mu}(B_R)}(\eps,\beta,T)\le\frac{18}5\eps\theta^2+14\,\eps\theta\sqrt{\frac{\mu}{L-\mu}} .
\end{equation*}
\item If $L_0\ge5\mu$ and $0<\eps\le\mu R^2/4$, then
\begin{equation*}
\Delta_{\cF_{\mathrm{SmSC}}^{L,\mu}(B_R)}(\eps,\beta,T)\le11\eps\theta^2+30\,\eps\theta\sqrt{\frac{\mu}{L-\mu}} .
\end{equation*}
\end{enumerate}
\end{theorem}

The bounds of Theorem~\ref{thm:main-smsc} do not depend on $R$. They are two explicit instances of the certificate of Lemma~\ref{lem:smsc-general} in Appendix~\ref{app:smsc}, which holds for all admissible parameters.

\begin{theorem}[Endpoint $L=\mu$]\label{thm:main-endpoint}
Let $\eps>0$ and $0<\beta<1/2$.
\begin{enumerate}[label=(\roman*),leftmargin=*]
\item The class $\cF_{\mathrm{SmSC}}^{\mu,\mu}(B_R)$ consists exactly of the functions $x\mapsto\frac\mu2\norm{x-x_c}^2+d$ on the ball $B_R$, where $x_c\in\R^n$ and $d\in\R$ are arbitrary.
\item A deterministic algorithm with $n+1$ queries on the sphere $\partial B_R$ returns, for every $\delta\ge0$, every $f\in\cF_{\mathrm{SmSC}}^{\mu,\mu}(B_R)$, and every $O\in\oo_\delta(f)$, a point $\widehat x\in B_R$ with
\[
f(\widehat x)-\min_{B_R}f\le\frac{n\delta^2}{2\mu R^2}.
\]
Hence $\Delta_{\cF_{\mathrm{SmSC}}^{\mu,\mu}(B_R)}(\eps,\beta,T)\ge R\sqrt{2\mu\eps/n}$ for every $T\ge n+1$.
\item If $T\le n$, then no randomized algorithm with $T$ queries $(\eps,\beta,0,T)$-solves $\cF_{\mathrm{SmSC}}^{\mu,\mu}(B_R)$; thus $\Delta_{\cF_{\mathrm{SmSC}}^{\mu,\mu}(B_R)}(\eps,\beta,T)=0$.
\item For every $L\ge\mu$, every $(n,T,\beta)$-hiding $\theta$, and $0<\eps\le\mu R^2/3$,
$\Delta_{\cF_{\mathrm{SmSC}}^{L,\mu}(B_R)}(\eps,\beta,T)\le R\sqrt{3\mu\eps}\,\theta$.
\end{enumerate}
In particular, for $T\ge n+1$, $\theta_*\le1/4$, and $0<\eps\le\mu R^2/3$,
\[
R\sqrt{\frac{2\mu\eps}{n}}\le\Delta_{\cF_{\mathrm{SmSC}}^{\mu,\mu}(B_R)}(\eps,\beta,T)\le R\sqrt{\frac{6\mu\eps\,\ell}{n-1}} .
\]
\end{theorem}

The endpoint MALN therefore grows linearly with $R$, while the bounds of Theorem~\ref{thm:main-smsc} are independent of $R$; the transition between the two regimes is quantified by Proposition~\ref{prop:near-end} and Fig.~\ref{fig:landscape}.

\subsection{Further Classes and the Neighborhoods of the Endpoints}

\begin{corollary}[H\"older continuous gradients]\label{cor:holder}
Let $\nu\in(0,1]$, $L_\nu>0$, $m_\nu:=(2^{\nu-1}L_\nu\eps^\nu)^{1/(1+\nu)}$, and assume $4\eps\le m_\nu R$. Then
\begin{equation*}
\Delta_{\cF_{\mathrm{H}}^{\nu,L_\nu}(B_R)}(\eps,\beta,T)
\le
\frac{12}5\eps\theta^2+7\cdot2^{\frac{1-\nu}{1+\nu}}\,\frac{\eps^{\frac{2+\nu}{1+\nu}}\,\theta}{R\,L_\nu^{\frac1{1+\nu}}},
\end{equation*}
and the same bound holds for $\cF_{\mathrm{H}}^{\nu,L_\nu}\cap\cF_{\mathrm{Lip}}^M$ whenever $M\ge m_\nu$.
\end{corollary}

For $\nu=1$ this is Theorem~\ref{thm:main-smooth}. As $\nu\downarrow0$, the exponent of $\eps$ tends to $2$, as in the Lipschitz case. This is an analogy between the scales, not an identification of the function classes at $\nu=0$: for a globally $M$-Lipschitz convex extension, differences of subgradients have norm at most $2M$, whereas the relative subdifferential at the boundary of $B_R$ contains unbounded normal components.

\begin{corollary}[Uniform convexity]\label{cor:uc}
Let $p\ge2$, $\mu_p>0$, $M_p:=2^{p-2}\mu_pR^{p-1}$, and $\rho_p:=(p\eps/(2^{p-1}\mu_p))^{1/p}$. If $\rho_p\le R$ and $M-M_p\ge4\eps/\rho_p$, then
\begin{equation*}
\Delta_{\cF_{\mathrm{UC}}^{M,\mu_p,p}(B_R)}(\eps,\beta,T)
\le
\frac{16}{7}\Bigl(\eps\theta^2+\frac{2\eps^2\theta}{(M-M_p)\rho_p}\Bigr).
\end{equation*}
\end{corollary}

For $p=2$ one has $M_2=\mu R$, $\rho_2=\sqrt{\eps/\mu}$, and $\cF_{\mathrm{UC}}^{M,\mu,2}(B_R)=\cF_{\mathrm{LipSC}}^{M,\mu}(B_R)$, so the corollary coincides with Theorem~\ref{thm:main-lipsc}. For $M\ge2M_p$ the term linear in $\theta$ is, up to factors depending on $p$, $\mu_p^{1/p}\eps^{2-1/p}\theta/M$, which for $\theta=\theta_*$ is the order of the level indicated for the $p$-growth condition in \cite{DvinskikhTominin2022}.

\begin{proposition}[Lipschitz strongly convex endpoint]\label{prop:lipsc-end}
Every $\mu$-strongly convex $M$-Lipschitz function on $B_R$ satisfies $M\ge\mu R$ and
\[
f(0)-\min_{B_R}f\le\frac{(M-\mu R)^2}{2\mu}.
\]
Consequently, if $M-\mu R\le\sqrt{2\mu\eps}$, the algorithm that makes no queries and outputs $0$ solves $\cF_{\mathrm{LipSC}}^{M,\mu}(B_R)$ at every noise level, and $\Delta_{\cF_{\mathrm{LipSC}}^{M,\mu}(B_R)}(\eps,\beta,T)=+\infty$.
\end{proposition}

Proposition~\ref{prop:lipsc-end} identifies a regime in which the class is solvable without queries. Theorem~\ref{thm:main-lipsc} gives a finite upper bound when $M-\mu R\ge4\sqrt{\mu\eps}$, using the remaining Lipschitz budget $M-\mu R$ in the construction. Together, these results demonstrate a transition near the endpoint, but they do not determine the optimal dependence on $M-\mu R$. The window $\sqrt{2\mu\eps}<M-\mu R<4\sqrt{\mu\eps}$ is covered by neither statement.

\begin{proposition}[Neighborhood of $L=\mu$]\label{prop:near-end}
Let $f\in\cF_{\mathrm{SmSC}}^{L,\mu}(B_R)$, $L_0:=L-\mu\ge0$, $O\in\oo_\delta(f)$, and $0<r\le R$. The deterministic algorithm that queries the $n+1$ vertices $ru_0,\ldots,ru_n$ of a regular simplex inscribed in $\partial B_r$, with $u_i$ as in Lemma~\ref{lem:simplex}, and outputs $\widehat x=\Proj_{B_R}(-\widehat c/\mu)$, where $\widehat c=\frac{n}{r(n+1)}\sum_iO(ru_i)u_i$, satisfies
\begin{equation}\label{eq:near-end}
f(\widehat x)-\min_{B_R}f\le\frac1{2\mu}\Bigl(\sqrt n\Bigl(\frac{L_0r}{4}+\frac\delta r\Bigr)+L_0R\Bigr)^2 .
\end{equation}
Consequently, for $T\ge n+1$ and $\Delta:=\Delta_{\cF_{\mathrm{SmSC}}^{L,\mu}(B_R)}(\eps,\beta,T)$,
\begin{align*}
&\Delta\ge R\sqrt{\frac{2\mu\eps}{n}}-L_0R^2\Bigl(\frac14+\frac1{\sqrt n}\Bigr),\\
&\Delta\ge\min\Bigl\{\frac{L_0R^2}{4},\ \frac{(\sqrt{2\mu\eps}-L_0R)^2}{nL_0}\Bigr\},
\end{align*}
the second inequality for $0<L_0R<\sqrt{2\mu\eps}$.
\end{proposition}

For example, if $T\ge n+1$ and $0\le L_0R\le\frac12\sqrt{\mu\eps/n}$, the first consequence gives $\Delta\ge\frac R2\sqrt{\mu\eps/n}$, since $\frac14+\frac1{\sqrt n}<1$ for $n\ge2$. Thus the endpoint scale persists up to absolute constants in this explicit neighborhood. A positive admissible noise level persists throughout $L_0R<\sqrt{2\mu\eps}$.

\subsection{Summary of Proved Scales}

In Table~\ref{tab:summary}, for $\theta=\theta_*$ the terms linear in $\theta$ carry $\theta\sqrt n=\Theta(\sqrt\ell)$ and the terms $\eps\theta^2$ carry $\ell$; for budgets polynomial in $n$, $\ell=O(\log n)$.

\begin{table*}[t]
\caption{Selected proved scales for the MALN. Absolute constants and, in the upper bounds, factors $\sqrt\ell$ on the terms of order $n^{-1/2}$ and $\ell$ on the terms $\eps/n$ are suppressed. Every row assumes $\eps>0$ and $0<\beta<1/2$. Every finite positive upper bound requires $n-1\ge32\ell$ for the actual budget $T$; this condition is not required for the zero or infinite cases. The lower-bound budgets $T_{\rm Lip}$, $T_{\rm SC}$, and $T_{\rm SmSC}$ are those of Propositions~\ref{prop:ach-lip}, \ref{prop:ach-sc}, and \ref{prop:ach-smsc}, respectively, with the parameters of the row. A dash means that no bound is asserted in that cell.}
\label{tab:summary}
\centering
\small
\setlength{\tabcolsep}{4pt}
\renewcommand{\arraystretch}{1.12}
\begin{tabular}{@{}>{\raggedright\arraybackslash}p{0.20\textwidth}>{\raggedright\arraybackslash}p{0.27\textwidth}>{\raggedright\arraybackslash}p{0.23\textwidth}>{\raggedright\arraybackslash}p{0.24\textwidth}@{}}
\toprule
Class & Joint parameter conditions & Upper scale & Lower scale and budget\\
\midrule
$\cF_{\mathrm{Lip}}^M$ & $\eps\le MR/8$ & $\frac{\eps^2}{\sqrt nMR}+\frac\eps n$ & $\frac{\eps^2}{\sqrt nMR}$; $T\ge T_{\rm Lip}$\\
\addlinespace
$\cF_{\mathrm H}^{\nu,L_\nu}\cap\cF_{\mathrm{Lip}}^{m_\nu}$ & $0<\nu\le1$, $4\eps\le m_\nu R$ & $\frac{\eps^2}{\sqrt n\,m_\nu R}+\frac\eps n$ & $\frac{\eps^2}{\sqrt n\,m_\nu R}$; $T\ge T_{\rm Lip}$, $M=m_\nu$\\
\addlinespace
$\cF_{\mathrm{Sm}}^L\cap\cF_{\mathrm{Lip}}^{\sqrt{L\eps}}$ & $\eps\le LR^2/16$ & $\frac{\eps^{3/2}}{\sqrt n\sqrt L R}+\frac\eps n$ & $\frac{\eps^{3/2}}{\sqrt n\sqrt L R}$; $T\ge T_{\rm Lip}$, $M=\sqrt{L\eps}$\\
\addlinespace
$\cF_{\mathrm{LipSC}}^{M,\mu}$ & $\eps\le\mu R^2$, $M\ge2\mu R$, $M-\mu R\ge4\sqrt{\mu\eps}$ & $\frac{\sqrt\mu\eps^{3/2}}{\sqrt nM}+\frac\eps n$ & $\frac{\sqrt\mu\eps^{3/2}}{\sqrt nM}$; $T\ge T_{\rm SC}$\\
\addlinespace
$\cF_{\mathrm{UC}}^{M,\mu_p,p}$ & $p\ge2$, $\rho_p\le R$, $M-M_p\ge4\eps/\rho_p$ & $\frac{\eps^2}{\sqrt n(M-M_p)\rho_p}+\frac\eps n$ & $\frac{\eps^2}{\sqrt nMR}$; $T\ge T_{\rm Lip}$\\
\addlinespace
$\cF_{\mathrm{SmSC}}^{L,\mu}$ & $L-\mu\ge5\mu$, $\eps\le\mu R^2/4$ & $\frac\eps{\sqrt n}\sqrt{\frac\mu{L-\mu}}+\frac\eps n$ & ---\\
\addlinespace
$\cF_{\mathrm{SmSC}}^{L,\mu,\circ}$ & $L-\mu\ge5\mu$, $\eps\le\mu R^2/6$ & $\frac\eps{\sqrt n}\sqrt{\frac\mu{L-\mu}}+\frac\eps n$ & $\frac\eps{\sqrt n}\sqrt{\frac\mu L}$; $T\ge T_{\rm SmSC}$\\
\addlinespace
$\cF_{\mathrm{SmSC}}^{L,\mu}$ & $0<(L-\mu)R\le\frac12\sqrt{\mu\eps/n}$, $\eps\le\mu R^2/3$ & $R\sqrt{\mu\eps/n}$ & $R\sqrt{\mu\eps/n}$; $T\ge n+1$\\
\addlinespace
$\cF_{\mathrm{SmSC}}^{\mu,\mu}$ & $\eps\le\mu R^2/3$ & $R\sqrt{\mu\eps/n}$ & $R\sqrt{\mu\eps/n}$; $T\ge n+1$\\
\addlinespace
$\cF_{\mathrm{SmSC}}^{\mu,\mu}$ & $T\le n$ & $0$ & $0$; exact-value unsolvability\\
\addlinespace
$\cF_{\mathrm{LipSC}}^{M,\mu}$ & $\mu R\le M\le\mu R+\sqrt{2\mu\eps}$ & $+\infty$ & $+\infty$; every $T\ge0$\\
\bottomrule
\end{tabular}
\par\vspace{5pt}
\parbox{0.97\textwidth}{\footnotesize The smooth and H\"older upper bounds also hold without the displayed Lipschitz restriction. Here $m_\nu$ is defined in Corollary~\ref{cor:holder}, and $M_p,\rho_p$ in Corollary~\ref{cor:uc}. The subclass $\cF_{\mathrm{SmSC}}^{L,\mu,\circ}$ is defined before Proposition~\ref{prop:ach-smsc}. Matching scales concern the terms of order $n^{-1/2}$ when those terms dominate. For $\cF_{\mathrm{SmSC}}^{L,\mu}$ without the restriction on the minimizer, the lower bound is open (Section~\ref{sec:discussion}). The lower bound for $\cF_{\mathrm{UC}}^{M,\mu_p,p}$ follows from the inclusion $\cF_{\mathrm{UC}}^{M,\mu_p,p}\subset\cF_{\mathrm{Lip}}^M$.}
\end{table*}

\section{Comparison With Achievable Noise Levels}\label{sec:matching}

Lower bounds on the MALN require algorithms. The following propositions are proved in Appendix~\ref{app:ach} in the model of Section~\ref{sec:model}, with queries in $B_R$ and an arbitrary fixed Borel perturbation.

\begin{proposition}[Two-point method, Lipschitz case]\label{prop:ach-lip}
Let $0<\eps\le MR$ and $0<\beta<1$. Let $J=\lceil\log_2(1/\beta)\rceil$ and $K=\lceil1200\,n^2M^2R^2/\eps^2\rceil$. For every $T\ge J(2K+1)$,
\[
\Delta_{\cF_{\mathrm{Lip}}^M(B_R)}(\eps,\beta,T)\ge\frac{\eps^2}{84\sqrt n\,MR}.
\]
\end{proposition}

\begin{proposition}[Two-point method, strongly convex case]\label{prop:ach-sc}
Let $0<\eps\le\mu R^2$ and $0<\beta<1$. Let $J=\lceil\log_2(1/\beta)\rceil$ and $K=\lceil217\,n^2M^2/(\mu\eps)\rceil$. For every $T\ge J(2K+1)$,
\[
\Delta_{\cF_{\mathrm{LipSC}}^{M,\mu}(B_R)}(\eps,\beta,T)\ge\frac{\sqrt\mu\,\eps^{3/2}}{19\sqrt n\,M}.
\]
\end{proposition}

Both methods use a symmetric two-point estimator along a random direction, project the iterates onto the ball $B_{R-h}$ with $h=7\eps/(64M)$ and $h=\eps/(8M)$, respectively, average them, and select the best of $J$ independent runs by the oracle values. The bias of the estimator caused by the perturbation has norm at most $\sqrt n\delta/h$; in the convex case this bias enters the optimization error through the diameter, and in the strongly convex case through its square divided by $\mu$.

For smooth objectives the smoothing bias is quadratic in $h$, which permits the radius $h\sim\sqrt{\eps/L}$ instead of $h\sim\eps/M$. The projection onto $B_{R-h}$ then costs nothing when the minimizer lies at distance at least $h$ from the boundary. Let $\cF_{\mathrm{SmSC}}^{L,\mu,\circ}(B_R)$ consist of the members of $\cF_{\mathrm{SmSC}}^{L,\mu}(B_R)$ whose minimizer over $B_R$ lies in $B_{R-\sqrt{\eps/L}}$; this subclass depends on $\eps$.

\begin{proposition}[Two-point method, smooth strongly convex case]\label{prop:ach-smsc}
Let $L\ge\mu$, $0<\eps\le\mu R^2$, and $0<\beta<1$. Let $J=\lceil\log_2(1/\beta)\rceil$ and $K=\lceil800\,n^2L^2R^2/(\mu\eps)\rceil$. For every $T\ge J(2K+1)$,
\[
\Delta_{\cF_{\mathrm{SmSC}}^{L,\mu,\circ}(B_R)}(\eps,\beta,T)\ge\frac{\eps}{5\sqrt n}\sqrt{\frac\mu L}.
\]
\end{proposition}

The order of Proposition~\ref{prop:ach-lip} is the level of \cite{DvinskikhTominin2022} and \cite{AlashqarGasnikovDvinskikhLobanov2023}, and the order of Proposition~\ref{prop:ach-sc} is the case $r=2$ of the level indicated in \cite{DvinskikhTominin2022}. The propositions state these orders in the model of Section~\ref{sec:model}, with queries restricted to $B_R$, an arbitrary fixed Borel perturbation, success probability $1-\beta$, and explicit constants, which Corollary~\ref{cor:two-sided} requires. Their purpose is to match noise tolerance; the elementary second-moment estimate used here yields budgets proportional to $n^2$, and these budgets are not claimed to be optimal in dimension.

\begin{corollary}[Two-sided bounds]\label{cor:two-sided}
Let $0<\beta<1/2$ and let $\theta$ be $(n,T,\beta)$-hiding with $T\ge J(2K+1)$ for the respective $J,K$ of Propositions~\ref{prop:ach-lip}--\ref{prop:ach-smsc}, and let $\kappa\ge1$.
\begin{enumerate}[label=(\alph*),leftmargin=*]
\item If $0<\eps\le MR/8$, then
\[
\frac{\eps^2}{84\sqrt nMR}\le\Delta_{\cF_{\mathrm{Lip}}^M(B_R)}(\eps,\beta,T)\le\frac43\eps\theta^2+\frac53\,\frac{\eps^2\theta}{MR}.
\]
\item If $M=\kappa\sqrt{L\eps}$ and $0<\eps\le LR^2/16$, then, with $\cF=\cF_{\mathrm{Sm}}^L(B_R)\cap\cF_{\mathrm{Lip}}^M(B_R)$,
\[
\frac{\eps^{3/2}}{84\kappa\sqrt n\sqrt LR}\le\Delta_{\cF}(\eps,\beta,T)\le\frac{12}5\eps\theta^2+\frac{7\eps^{3/2}\theta}{\sqrt LR}.
\]
\item If $0<\eps\le\mu R^2$, $M\ge2\mu R$, and $M-\mu R\ge4\sqrt{\mu\eps}$, then
\begin{align*}
\frac{\sqrt\mu\eps^{3/2}}{19\sqrt nM}&\le\Delta_{\cF_{\mathrm{LipSC}}^{M,\mu}(B_R)}(\eps,\beta,T)\\
&\le\frac{16}{7}\Bigl(\eps\theta^2+\frac{4\sqrt\mu\eps^{3/2}\theta}{M}\Bigr).
\end{align*}
\item If $\nu\in(0,1]$, $M=\kappa m_\nu$ with $m_\nu$ from Corollary~\ref{cor:holder}, and $4\eps\le m_\nu R$, then, with $\cF=\cF_{\mathrm{H}}^{\nu,L_\nu}(B_R)\cap\cF_{\mathrm{Lip}}^M(B_R)$,
\[
\frac{\eps^2}{84\kappa\sqrt n\,m_\nu R}\le\Delta_{\cF}(\eps,\beta,T)\le\frac{12}5\eps\theta^2+\frac{7\eps^2\theta}{m_\nu R}.
\]
\item If $L-\mu\ge5\mu$ and $0<\eps\le\mu R^2/6$, then, with $\cF=\cF_{\mathrm{SmSC}}^{L,\mu,\circ}(B_R)$,
\[
\frac{\eps}{5\sqrt n}\sqrt{\frac\mu L}\le\Delta_{\cF}(\eps,\beta,T)\le11\eps\theta^2+30\eps\theta\sqrt{\frac{\mu}{L-\mu}}.
\]
\end{enumerate}
\end{corollary}

Whenever the terms linear in $\theta$ dominate, the MALN is thus determined up to the factor $C\theta\sqrt n$ in (a), (c), and (e) and $C\kappa\theta\sqrt n$ in (b) and (d), where $C$ is the ratio of the coefficient of the term linear in $\theta$ in the upper bound to the coefficient of the lower bound: $C=\frac53\cdot84=140$ in (a), $C=7\cdot84=588$ in (b) and (d), $C=\frac{64}{7}\cdot19<174$ in (c), and $C=150\sqrt{L/(L-\mu)}\le150\sqrt{6/5}<165$ in (e). For $\theta=\theta_*$ one has $\theta\sqrt n=\sqrt{2\ell\,n/(n-1)}$. At the budget $T_0=J(2K+1)$, a hiding $\theta$ exists if $n-1\ge32\log(2(T_0+1)/(1-\beta))$ by \eqref{eq:theta-star}; for $MR/\eps=8$ and $\beta=1/4$ in (a) this holds for $n\ge870$. For (a), the optimality of the order $\eps^2/(\sqrt nMR)$ up to a logarithmic factor for $\eps\gtrsim MR/\sqrt n$ was observed in \cite{DvinskikhTominin2022} from the bound of \cite{LiRisteski2016}; Corollary~\ref{cor:two-sided} makes the constants explicit and gives the factor $\sqrt\ell$ when $\theta=\theta_*$ and $\kappa$ is bounded by an absolute constant, for every budget satisfying both the achievability and hiding conditions.

For $\theta=\theta_*$, the terms $\eps\theta_*^2=\Theta(\eps\ell/n)$ match, up to the factor $\ell$, the universal scale $\eps/n$ at which the simulated-annealing method of \cite{BelloniLiangNarayananRakhlin2015} optimizes functions uniformly close to convex ones with polynomially many queries under its stated geometric and regularity assumptions. This comparison does not supply an additional uniform lower bound, at the budget $T$ fixed here, for smooth classes with unrestricted affine parts. The smooth case in (b) is formulated on the subclass with Lipschitz constant proportional to $\sqrt{L\eps}$, which contains the hard functions of Theorem~\ref{thm:main-smooth}; on this subclass the Lipschitz method of Proposition~\ref{prop:ach-lip} already reaches the smooth scale, and (d) is the same statement for Corollary~\ref{cor:holder}. In (e) the subclass contains the hard functions of Theorem~\ref{thm:main-smsc}(b), whose minimizers lie at distance of order $\sqrt{\eps/\mu}$ from the origin; for minimizers on the boundary of $B_R$ the projection onto $B_{R-h}$ costs $\norm{\nabla f(x_f)}h$, which the smooth radius $h\sim\sqrt{\eps/L}$ does not control.

\begin{table}[t]
\caption{Noise models in related work and their role here.}
\label{tab:models}
\centering
\footnotesize
\setlength{\tabcolsep}{3pt}
\begin{tabular}{@{}>{\raggedright\arraybackslash}p{0.24\columnwidth}>{\raggedright\arraybackslash}p{0.4\columnwidth}>{\raggedright\arraybackslash}p{0.3\columnwidth}@{}}
\toprule
Source & Model & Role here\\
\midrule
\cite{SingerVondrak2015} & erroneous value oracle, convex & information-theoretic framework\\
\cite{LiRisteski2016} & fixed $F$ with $\sup|F-f|\le\delta$, accuracy measured on $F$, Lipschitz & same perturbations; Lipschitz scale\\
\cite{BelloniLiangNarayananRakhlin2015} & approximately convex $F$ & branch $\eps/n$\\
\cite{NesterovSpokoiny2017,Gasnikov2024,GorbunovDvurechenskyGasnikov2022} & bounded additive errors, often with stochastic components or unconstrained domains & algorithmic scales; smoothing\\
\cite{BayandinaGasnikovLagunovskaya2018,DvinskikhTominin2022,AlashqarGasnikovDvinskikhLobanov2023} & bounded deterministic noise, admissible level of two-point methods & orders of Propositions~\ref{prop:ach-lip}--\ref{prop:ach-sc}\\
\cite{LobanovBashirovGasnikov2024,GasnikovLobanovStonyakin2024} & higher-order smoothness or PL assumptions & comparison of smooth noise scales\\
\cite{Lobanov2023Stochastic} & stochastic errors independent of the random direction; bounded second moment & variance-based budget comparison\\
\cite{LobanovGasnikov2024LipschitzNoise} & Lipschitz bound on the perturbation & different noise parameter\\
\cite{DuchiJordanWainwrightWibisono2015,Shamir2017} & stochastic two-point feedback & estimator variance\\
\cite{AkhavanPontilTsybakov2021,AkhavanChzhenPontilTsybakov2024,BachPerchet2016} & adversarial or stochastic noise, higher-order smoothness & Section~\ref{sec:extensions}\\
\cite{Kerger2026,ZhangZhangQiLin2026} & exact values & budgets with vanishing MALN, Section~\ref{sec:matching}\\
This paper & fixed Borel perturbation, queries in $B_R$ & two-sided bounds\\
\bottomrule
\end{tabular}
\end{table}

The upper bounds also determine how many queries are needed at a given noise level.

\begin{corollary}[Query complexity under noise]\label{cor:budget}
Let $0<\eps\le MR/8$, $0<\beta<1/2$, $\delta>0$, and
\[
\theta(\delta):=\min\Bigl\{\frac14,\ \sqrt{\frac{3\delta}{8\eps}},\ \frac{3\delta MR}{10\eps^2}\Bigr\}.
\]
If a randomized algorithm with $T$ queries $(\eps,\beta,\delta,T)$-solves $\cF_{\mathrm{Lip}}^M(B_R)$, then
\[
T+1\ge\frac{1-\beta}{2}\exp\Bigl(\frac{(n-1)\,\theta(\delta)^2}{2}\Bigr).
\]
\end{corollary}

Since $\theta(\delta)^2=\min\{\frac1{16},\frac{3\delta}{8\eps},\frac{9\delta^2M^2R^2}{100\,\eps^4}\}$, the bound forces $T+1\ge\frac{1-\beta}2\exp(c\,n\min\{1,\delta/\eps,\delta^2M^2R^2/\eps^4\})$ with an absolute $c>0$ for $n\ge2$. By Proposition~\ref{prop:ach-lip}, polynomially many queries suffice for $\delta\le\eps^2/(84\sqrt nMR)$. For $\delta=t\,\eps^2/(\sqrt nMR)$ with $\delta\le\eps^3/(M^2R^2)$ the exponent is at least $c\,t^2$, so budgets polynomial in $n$ fail once $t\gg\sqrt{\log n}$; for $\delta\ge\eps$ the exponent is at least $c'n$ with an absolute $c'>0$, because $MR/\eps\ge8$.

Exact values are admissible for every $\delta\ge0$, so a budget for which exact values do not suffice is a budget at which the MALN vanishes. For randomized algorithms, \cite[Thm.~6]{ZhangZhangQiLin2026} applies to the success probability $1/2<1-\beta$ and to $2$-Lipschitz convex functions on the unit ball; after the rescaling $f\mapsto\frac{2}{MR}f(R\,\cdot)$, which maps $\cF_{\mathrm{Lip}}^M(B_R)$ onto this class and the accuracy $\eps$ to $\eps'=2\eps/(MR)$, it gives $\Delta_{\cF_{\mathrm{Lip}}^M(B_R)}(\eps,\beta,T)=0$ whenever $T<c\,n\min\{n,\eps'^{-2}\}/\log\min\{n,\eps'^{-2}\}$, $n\ge d_0$, and $\eps'\le\eps_0$, with the universal constants $c,d_0,\eps_0$ of that theorem. The deterministic lower bound \cite[Thm.~4]{Kerger2026} gives the same conclusion for the MALN over deterministic algorithms, that is, with $\Omega$ restricted to one-point spaces in Definition~\ref{def:maln}, at $\eps=\eps_0MR/\sqrt n$, $n\ge d_0$, and $T<c\,n^2/\log(n+1)$, with the universal constants $c,d_0,\eps_0$ of that theorem.

\section{Extensions and Limits of the Method}\label{sec:extensions}

\subsection{Effective Dimension}
The number of hidden directions, rather than the ambient dimension, governs the bounds.

\begin{corollary}[Functions of a $d$-dimensional projection]\label{cor:effdim}
Let $V\subset\R^n$ be a linear subspace of dimension $d\ge2$ with orthogonal projector $P$, and for $\cF\in\{\cF_{\mathrm{Lip}}^M,\cF_{\mathrm{Sm}}^L,\cF_{\mathrm{H}}^{\nu,L_\nu},\cF_{\mathrm{LipSC}}^{M,\mu},\cF_{\mathrm{SmSC}}^{L,\mu}\}$ let $\cF|_V$ consist of the members $f$ of $\cF(B_R)$ with $f-Q_\cF=(f-Q_\cF)\circ P$ on $B_R$, where $Q_\cF=Q$ for the strongly convex classes and $Q_\cF=0$ otherwise. Then Theorems~\ref{thm:main-lip}--\ref{thm:main-smsc} and Corollary~\ref{cor:holder} hold for $\cF|_V$ with $n$ replaced by $d$ in Definition~\ref{def:hiding}.
\end{corollary}

When $V$ is known, the lower bounds of Propositions~\ref{prop:ach-lip} and \ref{prop:ach-sc} also hold for $\cF|_V$ with $n$ replaced by $d$: every $f\in\cF|_V$ satisfies $f(Px)\le f(x)$ for $x\in B_R$, so its minimum over $B_R$ is attained on $V\cap B_R$, and the methods applied to the restriction of $f$ to the $d$-dimensional ball $V\cap B_R$ return points of $B_R$.

\subsection{Universal Convolution}
A transfer from the Lipschitz class to the smooth class through a translation-invariant convolution is limited by the following proposition, proved in Appendix~\ref{app:smooth}.

\begin{proposition}[Convolution barrier]\label{prop:barrier}
Let $U$ be a random vector in $\R^n$ with $\E\norm U<\infty$, $\gamma>0$, and $S_\gamma g(x):=\E g(x+\gamma U)$. If for every convex $M$-Lipschitz $g:\R^n\to\R$ the function $S_\gamma g$ has an $L$-Lipschitz gradient and $\sup|S_\gamma g-g|\le\mathsf b$, then $M\le2\sqrt{L\mathsf b}\,n^{-1/4}$.
\end{proposition}

Consequently, a Lipschitz hard family transferred to $\cF_{\mathrm{Sm}}^L$ by such a convolution with bias $\mathsf b=\Theta(\eps)$ has $M=O(\sqrt{L\eps}\,n^{-1/4})$, and the Lipschitz scale $\eps^2/(\sqrt nMR)$ is then at least of order $\eps^{3/2}/(n^{1/4}\sqrt LR)$. The Moreau envelope is a nonlinear operation whose bias $\lambda m^2/2$ is independent of the dimension. It does transfer information from the hidden cone to visible points, but only through the tangential component of its gradient, which Lemma~\ref{lem:moreau} bounds by the hidden slope $a$ rather than by the Lipschitz constant $m$. This is the single place where the smoothing step requires information specific to the family, and Fig.~\ref{fig:dimension} compares the two operations numerically.

\subsection{Higher-Order Smoothness}
Classes of convex functions with Lipschitz derivatives of order $p\ge2$ contain all linear functions, so the following elementary bound applies to them.

\begin{proposition}[Linear obstruction]\label{prop:linear}
Let $\cF$ be a class of functions on $B_R$ that contains $x\mapsto-\frac{3\eps}{2R}\ip wx$ for every $w\in\Sph$. For every $(n,T,\beta)$-hiding $\theta$, $\Delta_\cF(\eps,\beta,T)\le\frac32\eps\theta$.
\end{proposition}

\begin{proof}
Let $f_w(x)=-\frac{3\eps}{2R}\ip wx$ and let $O_w$ equal $0$ on $B_R\setminus A_w$ and $f_w$ on $A_w\cap B_R$. On visible points $|f_w|\le\frac32\eps\theta$, so $O_w\in\oo_{3\eps\theta/2}(f_w)$. The minimum of $f_w$ on $B_R$ is $-\frac32\eps$, while visible values are at least $-\frac32\eps\theta\ge-\frac38\eps$; hence every visible point has gap at least $\frac98\eps$, and every $\eps$-minimizer lies in $A_w$. Lemma~\ref{lem:hiding} with $H=0$ completes the proof.
\end{proof}

The Moreau envelope generally yields only $C^{1,1}$ regularity. Higher-order kernel methods \cite{BachPerchet2016,AkhavanPontilTsybakov2020,AkhavanChzhenPontilTsybakov2024} motivate studying smoother classes, but their guarantees depend on the oracle model, domain, and derivative assumptions. A heuristic balance $L_ph^{p+1}\sim\eps$ suggests $h=(\eps/L_p)^{1/(p+1)}$ and a candidate scale $\eps h/(\sqrt nR)$ when $h\le R$. This is not an upper or lower bound proved here. Any conjectured scale must also respect Proposition~\ref{prop:linear} and account for boundary queries and unrestricted lower-order derivatives. Constructing an appropriate hard family remains open.

\section{Numerical Illustrations}\label{sec:numerics}

The source code is available at \url{https://github.com/dmivilensky/MALN}. The envelope $e_\lambda g_w$ is evaluated through the dual formula of Lemma~\ref{lem:moreau}: by rotation invariance it depends on $x$ only through $(\ip wx,\norm{x-\ip wxw})$, the inner maximization is an explicit projection onto a planar section of $smK_w$, and the outer function of $s$ is concave and is maximized by golden-section search. Against direct numerical minimization of the primal problem in the invariant plane for 40 random parameter instances with 10 primal solver starts per instance, the maximal observed absolute difference is $1.3\cdot10^{-15}$.

\subsection{The Hiding Mechanism}
\begin{figure*}[t]
\centering
\includegraphics[width=0.92\textwidth]{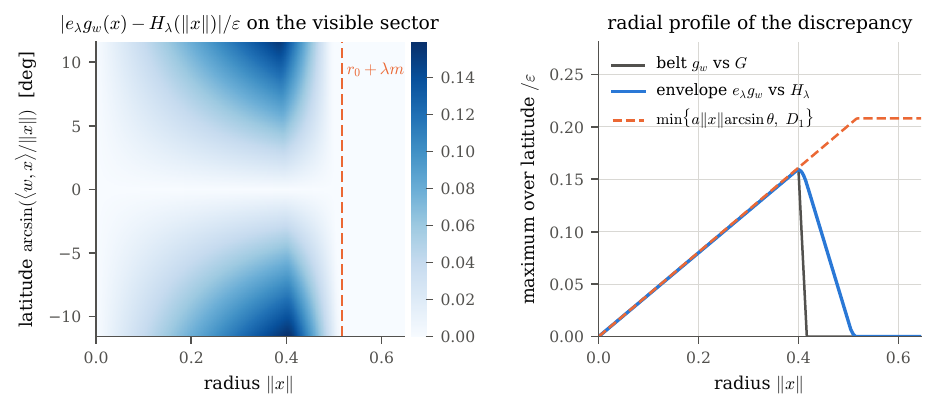}
\caption{Left: discrepancy $|e_\lambda g_w(x)-H_\lambda(\norm x)|/\eps$ on the visible sector as a function of the radius and the latitude $\arcsin(\ip wx/\norm x)\in[-\arcsin\theta,\arcsin\theta]$; the dashed line marks $r_0+\lambda m$. Right: its maximum over the latitude, the same quantity for the unsmoothed belt $|g_w-G|$, and the bound $\min\{a\norm x\arcsin\theta,D_1\}$. Parameters: $\eps=10^{-2}$, $L=R=1$, $\theta=0.2$.}
\label{fig:mechanism}
\end{figure*}

Fig.~\ref{fig:mechanism} shows the discrepancy between the envelope and its common radial branch on the visible sector for $\eps=10^{-2}$, $L=R=1$, and $\theta=0.2$ with the parameters of the proof of Theorem~\ref{thm:main-smooth}, that is, $m=0.1$, $a=0.02$, $r_0=0.417$, and $r_0+\lambda m=0.517$. The discrepancy vanishes for $\norm x\ge r_0+\lambda m$, where the computed values are below $3\cdot10^{-17}$. Its maximum $1.59\cdot10^{-3}$ is below the certified level $D_1=2.08\cdot10^{-3}$. The ratio of the discrepancy to the pointwise bound $a\norm x\arcsin(|\ip wx|/\norm x)$ of Steps 4--5 of the proof of Lemma~\ref{lem:moreau} reaches $1-5\cdot10^{-7}$ on the computational grid; at the maximizer of the discrepancy, the point of radius $0.400$ at the latitude $-\arcsin\theta$, it equals $0.987$, so the pointwise bound is nearly attained at these parameters. The discrepancy is asymmetric in the latitude. For the unsmoothed belt it vanishes beyond $b/m=0.400$ on the side of $+w$ and beyond $r_0=0.417$ on the side of $-w$, where the hidden branch $-a\ip wx$ is positive; the envelope extends both supports by $\lambda m=0.1$, and the maximum is attained on the side of $-w$ at the boundary latitude.

\subsection{Moreau Envelope Versus Universal Convolution}
\begin{figure}[t]
\centering
\includegraphics[width=\columnwidth]{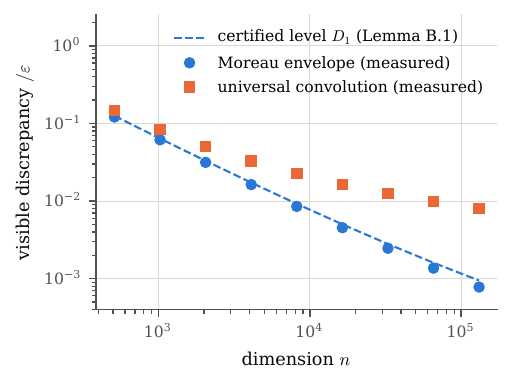}
\caption{Visible discrepancy of the smoothed belt family at the hiding parameter $\theta_*(n,10^6,1/4)$. Moreau envelope (circles) with its certified level $D_1$ (dashed), and universal convolution with the compact kernel $p_n$ at the same bias and certified smoothness (squares). Parameters: $\eps=3\cdot10^{-5}$, $L=R=1$.}
\label{fig:dimension}
\end{figure}

Fig.~\ref{fig:dimension} compares two smoothings of the same belt family at the same bias $\eps/2$ and the same certified gradient-Lipschitz constant $L$, for $\eps=3\cdot10^{-5}$, $L=R=1$, $T=10^6$, $\beta=1/4$, and $\theta=\theta_*(n,T,\beta)$. Both families have the hidden slope $a=2\eps/R$ of the proof of Theorem~\ref{thm:main-smooth}. For the Moreau envelope, $m=\sqrt{L\eps}$ and $\lambda=1/L$. For the convolution with the compact kernel $p_n(u)\propto(1-\norm u^2)^{\lceil n/2\rceil}\mathbf1\{\norm u<1\}$, the certified constant is $I_nm/\gamma$ with $I_n=\sup_{\norm v=1}\int|\partial_vp_n|$, which is computed in closed form and lies between $1.1283\sqrt n$ and $1.1287\sqrt n$ for the dimensions of the figure; matching the bias and the smoothness forces $m=\sqrt{L\eps/(2I_n\E\norm U)}$, in accordance with Proposition~\ref{prop:barrier}. The visible discrepancy of the convolution is computed by Monte Carlo integration with $4\cdot10^5$ common samples of $(U_1,U_2,\|U_{3:n}\|_2)$, the three scalar quantities entering the integrand, on a grid of $70$ radii and $11$ latitudes; the estimated pointwise standard error at each selected grid maximizer is below $10^{-5}$ of the reported value. This estimate does not include grid discretization error or the effect of selecting the maximum. The measured discrepancy of the envelope is at most $0.967\,D_1$ for all dimensions of the figure and at least $0.82\,D_1$. The ratio of the two discrepancies grows from $1.21$ at $n=2^9$ to $10.3$ at $n=2^{17}$, and over the last four dimensions the fitted slopes are $-0.85$ for the envelope and $-0.34$ for the convolution. The envelope slope lies between the slopes $-1$ and $-1/2$ of its two terms $\eps\theta^2$ and $\eps^{3/2}\theta/(\sqrt LR)$ and is close to the slope $-0.80$ of $D_1$ over the same dimensions, and the convolution inherits the flatter dependence induced by $m\sim n^{-1/4}$: its slope coincides with the slope $-0.34$ of its unsmoothed level $D_0$.

\subsection{The Smooth Strongly Convex Bracket}
\begin{figure}[t]
\centering
\includegraphics[width=\columnwidth]{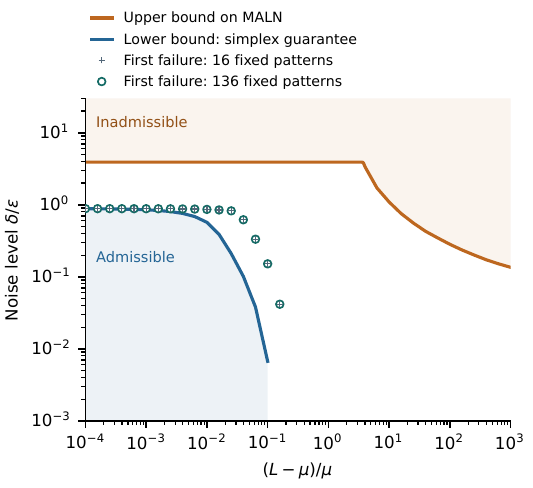}
\caption{Smooth strongly convex noise tolerance. Curves: the impossibility upper bound and the simplex lower guarantee. Crosses and circles: empirical first-failure levels on the same $40$ frozen quadratic objectives with $16$ and $136$ fixed perturbation patterns per objective, respectively; overlapping markers indicate agreement. Parameters: $n=256$, $\mu=R=1$, $\eps=10^{-2}$, $T=n+1$, $\beta=1/4$.}
\label{fig:landscape}
\end{figure}

Fig.~\ref{fig:landscape} displays the proved bounds on $\Delta_{\cF_{\mathrm{SmSC}}^{L,\mu}(B_R)}(\eps,\beta,T)$ as functions of $(L-\mu)/\mu$ for $n=256$, $\mu=R=1$, $\eps=10^{-2}$, $T=n+1$, $\beta=1/4$, and $\theta=\theta_*=0.226$. The upper bound is the minimum of Theorem~\ref{thm:main-endpoint}(iv) and the certificate of Lemma~\ref{lem:smsc-general} with parameters $(a,m,\rho)$ selected on a grid among those satisfying \eqref{eq:smsc-sep}; each sampled ordinate is obtained by evaluating an explicit bound. For Fig.~\ref{fig:landscape}, the selected parameters and strict separation margins are recorded in the repository and rechecked with high precision; plotted segments only interpolate between the evaluated points. The certificate falls below the endpoint upper bound between $(L-\mu)/\mu=3.5$ and $3.75$ and then decreases with $L-\mu$, from $3.4\eps$ at $(L-\mu)/\mu=4$ to $0.14\eps$ at $(L-\mu)/\mu=10^3$. The lower bound is Proposition~\ref{prop:near-end} with the query radius chosen optimally. At $L-\mu=10^{-4}\mu$ it equals $0.881\eps$, close to the endpoint value $R\sqrt{2\mu\eps/n}=0.884\eps$, while the endpoint upper bound is $3.92\eps$.

The markers report empirical first-failure levels on a frozen finite test family. We use $40$ quadratic objectives
\[
f(x)=Q(x)+\ip cx+\frac{L_0}{2}\norm{W^\top(x-p)}^2,
\]
where $W$ has orthonormal columns and $c$, $p$, and $W$ are drawn from the distribution specified in the repository. For each objective and each $L_0$, all perturbation patterns are fixed before varying $\delta$. The baseline test uses the signs of the simplex projections of four structured and four random directions, with both signs, giving $16$ patterns. An expanded test uses $64$ random directions instead, giving $136$ patterns; its first eight directions include the baseline test. For a pattern $z\in\{-1,1\}^{n+1}$ and $\delta>0$, the queried values are $f(ru_i)+\delta z_i$, with $r=\min\{R,2\sqrt{\delta/L_0}\}$. These values extend to a fixed Borel oracle in $\oo_\delta(f)$.

For each frozen objective and pattern we locate the first crossing of the loss level $\eps$ as $\delta$ increases, and take the smallest such crossing over the test family. We do not assume that the loss is monotone. On either side of $r=R$, the reconstructed linear coefficient is affine in $\sqrt\delta$ or in $\delta$; after splitting the projected output into its interior and boundary regimes, the level-crossing candidates are roots of quadratic or quartic polynomials. Candidates introduced by squaring are discarded using the original loss, and crossing roots are refined in local brackets. The code in the repository independently checks the losses immediately before and after $3915$ selected crossings and at smaller noise levels, and stores a failing instance for every plotted empirical point. Zero first-failure levels are retained in the data but omitted from the logarithmic plot.

At $L_0=10^{-4}\mu$, both tests give $0.883712\eps$, within $0.02\%$ of the endpoint guarantee $R\sqrt{2\mu\eps/n}=0.883883\eps$. At $L_0=0.1\mu$, the first-failure level divided by the guarantee of Proposition~\ref{prop:near-end} is $22.5818$ in both tests. These numbers describe this fixed test family and radius rule; they are neither uniform guarantees over all bounded perturbations nor the true worst-case tolerance of the simplex method. The theoretical curves are computed independently of these markers.

\subsection{Breakdown of the Two-Point Method}\label{sec:breakdown}
The guaranteed level of Proposition~\ref{prop:ach-lip} can be compared with the noise level at which the method breaks down on an explicit instance. Let $u\in\Sph$ and $f(x)=-\frac{2\eps}R\ip ux$, which belongs to $\cF_{\mathrm{Lip}}^M(B_R)$ for $\eps\le MR/2$, and let $O(x)=f(x)+\delta\operatorname{sgn}\ip ux$, where $\operatorname{sgn}(t)=1$ for $t>0$ and $\operatorname{sgn}(t)=-1$ for $t\le0$; then $O\in\oo_\delta(f)$. For the estimator $\widehat G$ of Appendix~\ref{app:ach}, with $e$ uniform on $\Sph$, $x\in B_{R'}$, $R'=R-h$, $\ip ux=s$, and $\tau_u=\ip ue$,
\[
\E\ip{\widehat G(x,e)}u=-\frac{2\eps}R+\frac{n\delta}h\,\E\bigl[|\tau_u|\,\mathbf 1\{|\tau_u|>|s|/h\}\bigr],
\]
because $\operatorname{sgn}(s+h\tau_u)-\operatorname{sgn}(s-h\tau_u)=2\operatorname{sgn}(\tau_u)\mathbf 1\{|\tau_u|>|s|/h\}$ almost surely and $\E \tau_u^2=1/n$. At $s=0$, the conditional mean increment of $\ip u{x_k}$ before projection changes sign at
\begin{gather*}
\delta_{\mathrm{b}}:=\frac{2\eps h}{nR\,\E|\tau_u|}=\frac{7}{32}\,\frac{\eps^2}{MR\,n\E|\tau_u|},\\
\E|\tau_u|=\frac{\Gamma(n/2)}{\sqrt\pi\,\Gamma((n+1)/2)},
\end{gather*}
with $h=7\eps/(64M)$ from the proof of Proposition~\ref{prop:ach-lip}. Since $\sqrt n\,\E|\tau_u|\to\sqrt{2/\pi}$, $\delta_{\mathrm{b}}\sqrt nMR/\eps^2\to\frac7{32}\sqrt{\pi/2}=0.274$, and the ratio of $\delta_{\mathrm{b}}$ to the guaranteed level $\eps^2/(84\sqrt nMR)$ equals $\frac{147}8(\sqrt n\,\E|\tau_u|)^{-1}$: it is $20.4$ for $n=2$, tends to $\frac{147}8\sqrt{\pi/2}$, and is smaller than $23.1$ for every $n$, because $\Gamma(x+\frac12)\le\sqrt x\,\Gamma(x)$ for $x>0$ by the log-convexity of $\Gamma$. For this instance the iteration of Lemma~\ref{lem:sgd}(a) depends on $x_k$ only through $\ip u{x_k}$ and $\norm{x_k-\ip u{x_k}u}$, and on $e_k$ only through $\ip u{e_k}$ and the component of $e_k$ along $x_k-\ip u{x_k}u$, whose joint law is that of two coordinates of a uniform point of $\Sph$; the code in the repository simulates this chain exactly in distribution, with all parameters of Proposition~\ref{prop:ach-lip}, and compares it with the full $n$-dimensional method for $n=4$. Table~\ref{tab:breakdown} reports single runs of $2K$ queries for $M=R=1$ and $\eps=1/8$ at $\delta=0.98\,\delta_{\mathrm{b}}$ and $\delta=1.02\,\delta_{\mathrm{b}}$. In every run at $0.98\,\delta_{\mathrm{b}}$ the returned point has gap at most $0.049\eps$, and in every run at $1.02\,\delta_{\mathrm{b}}$ the averaged iterate remains at the discontinuity and has gap larger than $2\eps$. For these seeds and iteration budgets, the observed transition between successful and unsuccessful runs lies between $0.98\,\delta_{\mathrm{b}}$ and $1.02\,\delta_{\mathrm{b}}$. These single-run observations neither establish a worst-case failure probability for the repeated-run selection procedure nor prove an impossibility result for other algorithms.

\begin{table}[t]
\caption{Two-point method of Proposition~\ref{prop:ach-lip} on the instance of Section~\ref{sec:breakdown}: iterations per run, breakdown level, and gaps of the points returned by single runs (maximum over runs at $0.98\,\delta_{\mathrm{b}}$, minimum over runs at $1.02\,\delta_{\mathrm{b}}$; four runs for $n\le32$, two for $n=64$). Parameters: $M=R=1$, $\eps=1/8$.}
\label{tab:breakdown}
\centering
\footnotesize
\setlength{\tabcolsep}{2.4pt}
\begin{tabular}{@{}rrcccc@{}}
\toprule
$n$ & $K$ & $\frac{\delta_{\mathrm{b}}\sqrt nMR}{\eps^2}$ & $\frac{84\sqrt nMR\,\delta_{\mathrm{b}}}{\eps^2}$ & gap$/\eps$, $0.98\delta_{\mathrm{b}}$ & gap$/\eps$, $1.02\delta_{\mathrm{b}}$\\
\midrule
$2$ & $3.07\cdot10^{5}$ & $0.243$ & $20.41$ & $0.048$ & $2.005$\\
$4$ & $1.23\cdot10^{6}$ & $0.258$ & $21.65$ & $0.049$ & $2.003$\\
$8$ & $4.92\cdot10^{6}$ & $0.266$ & $22.32$ & $0.046$ & $2.002$\\
$16$ & $1.97\cdot10^{7}$ & $0.270$ & $22.67$ & $0.044$ & $2.001$\\
$32$ & $7.86\cdot10^{7}$ & $0.272$ & $22.85$ & $0.044$ & $2.001$\\
$64$ & $3.15\cdot10^{8}$ & $0.273$ & $22.94$ & $0.043$ & $2.001$\\
\bottomrule
\end{tabular}
\end{table}

\section{Discussion}\label{sec:discussion}

A single explicit family yields the upper bounds for the Lipschitz, strongly convex, uniformly convex, smooth, H\"older, and smooth strongly convex classes. On visible points beyond an explicit radius, the belt family agrees exactly with a common radial function, and the Moreau envelope preserves this property with a discrepancy controlled by the hidden slope. The two-sided results, with $\theta=\theta_*$ and bounded proportionality constants for the Lipschitz subclasses, have the following scope: their terms of order $n^{-1/2}$ coincide with the noise tolerance of the two-point methods of Section~\ref{sec:matching} up to the factor $\sqrt\ell$ and absolute constants on the Lipschitz class, on the Lipschitz strongly convex class for $M\ge2\mu R$ and $M-\mu R\ge4\sqrt{\mu\eps}$, on the smooth and H\"older classes with Lipschitz constants proportional to that of the hard family, and on the smooth strongly convex class for $L-\mu\ge5\mu$ and minimizers at distance at least $\sqrt{\eps/L}$ from the boundary; these terms dominate for $\eps/(MR)\gtrsim\sqrt{\ell/n}$ in the Lipschitz case, and the terms $\eps\theta^2$ correspond, up to $\ell$, to the algorithmic scale $\eps/n$ of \cite{BelloniLiangNarayananRakhlin2015}. For $T\ge n+1$, $n-1\ge32\ell$, and $\eps\le\mu R^2/3$ the endpoint $L=\mu$ is determined up to $O(\sqrt\ell)$, and the query threshold $n+1$ is exact for randomized algorithms.

For the smooth strongly convex class, Theorem~\ref{thm:main-endpoint} gives two-sided bounds at $L=\mu$ in its stated budget, dimension, and accuracy ranges. Corollary~\ref{cor:two-sided}(e) gives matching leading terms on $\cF_{\mathrm{SmSC}}^{L,\mu,\circ}(B_R)$ for $L-\mu\ge5\mu$ and $\eps\le\mu R^2/6$, provided the terms linear in $\theta$ dominate. Between these regimes, Lemma~\ref{lem:smsc-general} supplies additional parameter-dependent upper bounds, including the numerical improvements in Fig.~\ref{fig:landscape}, but a sharp explicit characterization is not established. The lower bounds available between these regimes are those of Proposition~\ref{prop:near-end}, which decrease as $L-\mu$ grows, and of Proposition~\ref{prop:ach-smsc} on $\cF_{\mathrm{SmSC}}^{L,\mu,\circ}(B_R)$. For minimizers near the boundary, the contraction used by the two-point method incurs a first-order error controlled by $\norm{\nabla f(x_f)}h$, which is not uniformly bounded by smoothness alone.

A candidate hard family for $L-\mu\ll\mu$ is the envelope $e_{1/L}(Q_L+g_w)$ with $Q_L(x)=\frac{\mu L}{2(L-\mu)}\norm x^2$: the envelope with parameter $1/L$ of a $\frac{\mu L}{L-\mu}$-strongly convex function is $L$-smooth and $\mu$-strongly convex for every $L>\mu$. Further directions are the higher-order classes of Section~\ref{sec:extensions}, the exact power of the logarithm, and matching algorithms for the smooth classes without a Lipschitz bound proportional to $\sqrt{L\eps}$; kernel-based methods reach the order $\eps^{3/2}/\sqrt n$ for objectives of higher-order smoothness under bounded-variance noise in the unconstrained setting \cite{LobanovBashirovGasnikov2024}. The exact-value lower bounds of \cite{Kerger2026,ZhangZhangQiLin2026} rely on max-affine and random support-function families whose individual pieces carry small signals; quantifying how deterministic perturbations of size $\delta$ enlarge the corresponding query complexity would connect the polynomial and the exponential regimes of Corollary~\ref{cor:budget}.

\appendices
\section*{Guide to the Appendices}
Appendix~\ref{app:lip} proves the hiding lemma and the Lipschitz, strongly convex, and uniformly convex bounds.
Appendix~\ref{app:smooth} establishes Moreau hiding and the smooth and H\"older bounds.
Appendix~\ref{app:smsc} gives the smooth strongly convex certificate.
Appendix~\ref{app:endpoint} treats quadratic rigidity, the exact query threshold, and the neighborhood of $L=\mu$.
Appendix~\ref{app:ach} proves achievability, the budget bound, and the effective-dimension reduction.

\section{Hiding Lemma, Belt Family, and Lipschitz Classes}\label{app:lip}

\begin{lemma}[Spherical caps]\label{lem:cap}
Let $w$ be uniform on $\Sph$, $\theta\in(0,1)$, and $x\ne0$. Then $\Prob\{|\ip wx|>\theta\norm x\}\le2(1-\theta^2)^{(n-1)/2}$.
\end{lemma}

\begin{proof}
By rotation invariance we may take $x=e_1$. Write $w=Z/\norm Z$ with $Z$ standard Gaussian in $\R^n$ and $S_{n-1}:=\sum_{i\ge2}Z_i^2$, a $\chi^2_{n-1}$ variable independent of $Z_1$. The event $|w_1|>\theta$ is $Z_1^2>cS_{n-1}$ with $c=\theta^2/(1-\theta^2)$. Since $\Prob\{|Z_1|>u\}\le2e^{-u^2/2}$ for $u\ge0$,
\begin{align*}
\Prob\{|w_1|>\theta\}&\le2\,\E e^{-cS_{n-1}/2}=2(1+c)^{-(n-1)/2}\\
&=2(1-\theta^2)^{(n-1)/2}.\qedhere
\end{align*}
\end{proof}

\begin{proof}[Proof of Lemma~\ref{lem:hiding}]
Fix an algorithm $(\Omega,\Prob_A,(\chi_t))$ with $T$ queries. For $\omega\in\Omega$ define the reference points $\xi_t(\omega):=\chi_t(\omega,H(\xi_1(\omega)),\ldots,H(\xi_{t-1}(\omega)))$ for $t=1,\ldots,T+1$; they are Borel in $\omega$ and do not depend on $w$. Fix $w$ and $\omega$ with $\xi_t(\omega)\notin E_w$ for all $t\le T+1$. By induction on $t$, the queries of the algorithm against $O_w$ are $x_t=\xi_t(\omega)$: if $x_s=\xi_s(\omega)$ for $s<t$, then $y_s=O_w(\xi_s(\omega))=H(\xi_s(\omega))$ by (b), so $x_t=\xi_t(\omega)$. Hence the output is $\xi_{T+1}(\omega)\notin E_w$, and by (c) it is not an $\eps$-minimizer of $f_w$. For each fixed $w$ the success event is Borel in $\omega$, because $f_w$, $O_w$, and the maps $\chi_t$ are Borel, and it is contained in $\bigcup_{t\le T+1}\{\omega:\xi_t(\omega)\in E_w\}$. Therefore
\[
\Prob_A\{\text{success on }(f_w,O_w)\}\le\sum_{t=1}^{T+1}\Prob_A\{\xi_t\in E_w\}.
\]
The right side is Borel in $w$ by the joint measurability of $\mathbf 1\{x\in E_w\}$, and by Fubini's theorem and (d)
\begin{align*}
&\int\sum_{t=1}^{T+1}\Prob_A\{\xi_t\in E_w\}\,d\Prob_W(w)\\
&\quad=\sum_{t=1}^{T+1}\int\Prob_W\{\xi_t(\omega)\in E_w\}\,d\Prob_A(\omega)\le(T+1)p .
\end{align*}
If the algorithm $(\eps,\beta,D,T)$-solved $\cF$, then by (a) the integrand on the left side of the previous display would be at least $1-\beta$ for every $w$, hence so would the integral, which contradicts $(T+1)p<1-\beta$. The bound on $\Delta_\cF$ follows from the monotonicity noted after the definition of MALN. For the double cones, $0\notin A_w$, $(w,x)\mapsto\mathbf 1\{x\in A_w\}$ is the indicator of an open set, and Lemma~\ref{lem:cap} gives (d) with $p=\min\{1,2(1-\theta^2)^{(n-1)/2}\}$.
\end{proof}

\begin{lemma}[Support function of the belt]\label{lem:support}
For $x\in\R^n$ let $r=\norm x$, $t=\ip wx$, and $z=\norm{x-tw}$. Then
\[
\sigma_w(x)=
\begin{cases}
r,& |t|\le\theta r,\\
\theta|t|+\sqrt{1-\theta^2}\,z,& |t|>\theta r .
\end{cases}
\]
Moreover $\sigma_w$ is convex and $1$-Lipschitz, $\theta r\le\sigma_w(x)\le r$, $\sigma_{\Theta w}(\Theta x)=\sigma_w(x)$ for every orthogonal $\Theta$, and for every visible $x\ne0$ one has $x/\norm x\in K_w$ and $x/\norm x\in\partial\sigma_w(x)$.
\end{lemma}

\begin{proof}
Write $v\in K_w$ as $v=sw+v_\perp$ with $v_\perp\perp w$, $|s|\le\theta$, and $\norm{v_\perp}^2\le1-s^2$. Then $\ip vx=st+\ip{v_\perp}{x-tw}\le st+\sqrt{1-s^2}\,z$, with equality for a suitable $v_\perp$. The function $\varphi(s)=st+\sqrt{1-s^2}z$ is concave on $[-1,1]$, and by the Cauchy--Schwarz inequality $\varphi(s)\le r$ with equality at $s=t/r$ when $r>0$. If $|t|\le\theta r$ this maximizer is feasible and $\sigma_w(x)=r$. Otherwise concavity places the maximum over $[-\theta,\theta]$ at $s=\theta\operatorname{sgn}t$, which gives the second formula. As a support function of a subset of the unit ball, $\sigma_w$ is convex, $1$-Lipschitz, and at most $r$; the choice $v=\theta x/r\in K_w$ gives $\sigma_w(x)\ge\theta r$. Invariance follows from $\Theta K_w=K_{\Theta w}$. For visible $x\ne0$ the vector $x/\norm x$ lies in $K_w$ and attains $\ip vx=r=\sigma_w(x)$; the maximizers of $v\mapsto\ip vx$ over $K_w$ are subgradients of $\sigma_w$ at $x$.
\end{proof}

\begin{lemma}[Belt family]\label{lem:belt}
Let $0<a\le m$, $\rho>0$, $\theta\in(0,1/4]$, and let $g_w,G,r_0,D_0$ be defined by \eqref{eq:belt}--\eqref{eq:D0}. Then:
\begin{enumerate}[label=(\roman*),leftmargin=*]
\item $g_w$ is convex and $m$-Lipschitz on $\R^n$, and $g_{\Theta w}(\Theta x)=g_w(x)$ for every orthogonal $\Theta$;
\item $\min_{\R^n}g_w=g_w(\rho w)=-a\rho$;
\item every visible $x$ satisfies $|g_w(x)-G(x)|\le a\theta\min\{\norm x,r_0\}\le D_0$, and if $\norm x\ge r_0$, then $g_w(x)=G(x)=m\norm x-b$.
\end{enumerate}
\end{lemma}

\begin{proof}
(i) $g_w$ is the maximum of the convex $m$-Lipschitz function $m\sigma_w-b$ and the linear $a$-Lipschitz function $-a\ip w\cdot$, with $a\le m$; invariance follows from Lemma~\ref{lem:support}. (ii) For $\norm x=r$, Lemma~\ref{lem:support} gives $g_w(x)\ge\max\{m\theta r-b,-ar\}$. The first expression increases and the second decreases in $r$, and they coincide at $r=b/(m\theta+a)=\rho$, where both equal $-a\rho$. At $x=\rho w$ one has $\sigma_w(\rho w)=\theta\rho$, so both branches equal $-a\rho$. (iii) For visible $x$, $\sigma_w(x)=\norm x$ and $|a\ip wx|\le a\theta\norm x$. Since $u\mapsto\max\{c,u\}$ is $1$-Lipschitz, $|g_w(x)-G(x)|=|\max\{m\norm x-b,-a\ip wx\}-\max\{m\norm x-b,0\}|\le a\theta\norm x$. If $\norm x\ge r_0$, then $m\norm x-b\ge a\theta\norm x\ge\max\{0,-a\ip wx\}$, so both maxima equal $m\norm x-b$. Finally $m>a\theta$ because $a\le m$ and $\theta<1$, so $r_0$ is finite.
\end{proof}

\begin{proof}[Proof of Theorem~\ref{thm:main-lip}]
Take $m=M$, $\rho=R$, $a=5\eps/(4R)$; then $a\le m$ and $q=a/m=5\eps/(4MR)\le\frac5{32}$. Let $f_w=g_w|_{B_R}\in\cF_{\mathrm{Lip}}^M(B_R)$, let $O_w$ equal $G$ on visible points of $B_R$ and $g_w$ on $A_w\cap B_R$, and let $H=G$. By Lemma~\ref{lem:belt}, $O_w\in\oo_{D_0}(f_w)$ and $\min_{B_R}f_w=f_w(Rw)=-\frac54\eps$. By \eqref{eq:r0-ratio}, $D_0=a\rho\,\theta\frac{\theta+q}{1-q\theta}=\frac54\eps\theta\frac{\theta+q}{1-q\theta}$, and $\theta\le\frac14$, $q\le\frac5{32}$ give $1-q\theta\ge15/16$. Therefore $D_0\le\frac43\eps\theta(\theta+q)\le\eps/6$, where the last step uses $q\le1/4$ and $\theta\le1/4$. Every visible $x$ satisfies $f_w(x)\ge G(x)-D_0\ge-\eps/6$, so its gap is at least $(5/4-1/6)\eps>\eps$, which is (c) of Lemma~\ref{lem:hiding} with $E_w=A_w\cap B_R$. The functions $f_w$ and $O_w$ are Borel. Lemma~\ref{lem:hiding} shows that no algorithm solves $\cF_{\mathrm{Lip}}^M(B_R)$ at level $D_0$, hence at any level $D\ge D_0$, and $\frac43\eps\theta(\theta+q)=\frac43\eps\theta^2+\frac53\eps^2\theta/(MR)$ is the stated bound. Substituting $\theta_*^2=2\ell/(n-1)$ gives the last formula.
\end{proof}

\begin{proof}[Proof of Theorem~\ref{thm:main-lipsc}]
Take $m=M-\mu R$, $\rho=\sqrt{\eps/\mu}\le R$, and $a=2\sqrt{\mu\eps}$, so that $a\rho=2\eps$ and $q=a/m\le\frac12$. Let $f_w=(g_w+Q)|_{B_R}$. Then $f_w-Q=g_w$ is convex, and $f_w$ is $(m+\mu R)=M$-Lipschitz on $B_R$, so $f_w\in\cF_{\mathrm{LipSC}}^{M,\mu}(B_R)$. Let $O_w$ equal $G+Q$ on visible points and $f_w$ on $A_w\cap B_R$, and $H=G+Q$; then $O_w\in\oo_{D_0}(f_w)$. By \eqref{eq:r0-ratio} and $1-q\theta\ge\frac78$, $D_0=2\eps\theta\frac{\theta+q}{1-q\theta}\le\frac{16}7\eps\theta(\theta+q)\le\frac{16}7\cdot\frac14\cdot\frac34\eps=\frac37\eps$. Now $\min_{B_R}f_w\le f_w(\rho w)=-2\eps+\frac\eps2$, and every visible $x$ satisfies $f_w(x)\ge g_w(x)\ge-D_0\ge-\frac37\eps$, so its gap is strictly larger than $\eps$. Lemma~\ref{lem:hiding} applies with $D=D_0\le\frac{16}7\eps\theta(\theta+q)$ and $q=2\sqrt{\mu\eps}/(M-\mu R)$.
\end{proof}

\begin{proof}[Proof of Corollary~\ref{cor:uc}]
Let $\psi_p(x):=\frac{2^{p-2}\mu_p}{p}\norm x^p$. By \cite[Lemma~4]{Nesterov2008Cubic}, $\frac1p\norm y^p\ge\frac1p\norm x^p+\ip{\norm x^{p-2}x}{y-x}+\frac1p2^{2-p}\norm{y-x}^p$, hence $\psi_p(y)\ge\psi_p(x)+\ip{\nabla\psi_p(x)}{y-x}+\frac{\mu_p}p\norm{y-x}^p$. Take $m=M-M_p$, $\rho=\rho_p$, $a=2\eps/\rho_p\le m/2$, and $f_w=(g_w+\psi_p)|_{B_R}$. On $B_R$, $\norm{\nabla\psi_p}\le M_p$, so $f_w$ is $M$-Lipschitz. Since $g_w$ and $\psi_p$ are finite and convex on $\R^n$, the subdifferential of $f_w$ relative to $B_R$ at $x$ is $\partial g_w(x)+\nabla\psi_p(x)+N_{B_R}(x)$ \cite{BauschkeCombettes2017}, where $N_{B_R}(x)$ is the normal cone. For $s=s_g+\nabla\psi_p(x)+\varsigma$ in this set and $y\in B_R$, adding $g_w(y)\ge g_w(x)+\ip{s_g}{y-x}$, the inequality for $\psi_p$, and $\ip{\varsigma}{y-x}\le0$ gives the defining inequality of $\cF_{\mathrm{UC}}^{M,\mu_p,p}$. The value $\psi_p(\rho_pw)=\eps/2$ gives $\min f_w\le-\frac32\eps$, visible values are at least $-D_0\ge-\frac37\eps$ by the computation in the proof of Theorem~\ref{thm:main-lipsc}, and the oracle equal to $G+\psi_p$ on visible points has error $D_0\le\frac{16}7\eps\theta(\theta+q)$. Lemma~\ref{lem:hiding} gives the bound with $q=2\eps/(\rho_p(M-M_p))$. For $p=2$, suppose first that $f$ satisfies the uniform-convexity inequality. At every interior point $x$, choose $s\in\partial f(x)$ and subtract the identity for $Q$ to obtain
\[
(f-Q)(y)\ge(f-Q)(x)+\ip{s-\mu x}{y-x}.
\]
Thus $f-Q$ has an affine supporting function at each interior point and is convex there; continuity extends convexity to $B_R$. Conversely, if $f-Q$ is convex, the convex subdifferential sum rule for $(f-Q)+I_{B_R}$ and the smooth function $Q$ identifies every relative subgradient of $f$ as $s_h+\mu x$, where $s_h\in\partial((f-Q)+I_{B_R})(x)$. Adding the supporting inequality for $f-Q$ to the quadratic identity gives the defining inequality for $p=2$. Here $I_{B_R}$ is zero on $B_R$ and $+\infty$ outside.
\end{proof}

\begin{proof}[Proof of Proposition~\ref{prop:lipsc-end}]
Fix $u\in\Sph$ and let $\varphi(t)=f(tu)$ and $\psi(t)=\varphi(t)-\frac\mu2t^2$ on $[-R,R]$. The function $\psi$ is convex and $\varphi$ is $M$-Lipschitz, so for $-R<s<t<R$ the one-sided derivatives satisfy $\psi'_+(s)\le\psi'_-(t)=\varphi'_-(t)-\mu t\le M-\mu t$. Letting $t\uparrow R$ gives $\psi'_+(s)\le\varpi:=M-\mu R$, and symmetrically $\psi'_-(s)\ge-\varpi$. Since $\psi'_-\le\psi'_+$, this forces $\varpi\ge0$, and it shows that $\psi$ is $\varpi$-Lipschitz on $(-R,R)$, hence on $[-R,R]$. As $u$ is arbitrary, $\tilde f:=f-Q$ satisfies $\tilde f(x)\ge\tilde f(0)-\varpi\norm x$ on $B_R$, and
\[
f(x)-f(0)\ge\frac\mu2\norm x^2-\varpi\norm x\ge-\frac{\varpi^2}{2\mu}.
\]
If $\varpi\le\sqrt{2\mu\eps}$, the output $0$ is an $\eps$-minimizer of every member of the class, independently of the oracle.
\end{proof}

\begin{remark}[Relation to the construction of Li and Risteski]\label{rem:li-risteski}
The construction of \cite[Sec.~4]{LiRisteski2016} on the ball $B_{1/2}$ takes a hidden direction $w$ with $\norm w=\frac12$, the linear branch $x\mapsto-8\eps\ip wx$, and the maximum $f_w$ of this branch and the lower convex envelope of $\frac12\norm x^{1+\alpha}$ restricted to the part of the ball outside a core ball and outside a double cone around $w$; the perturbed function equals $\max\{f_w,\delta/2\}$ on the visible part of the core and $f_w$ elsewhere, and it coincides with a function independent of $w$ outside the double cone. The belt family shares the hidden linear branch in a random double cone and the common oracle outside it. Its radial part $m\sigma_w-b$ equals $m\norm x-b$ on the entire visible region, so Lemma~\ref{lem:belt}(iii) gives exact commonness beyond the explicit radius $r_0$ without a truncation at the level $\delta/2$, the discrepancy bound $D_0=a\theta r_0$, and a global definition on $\R^n$. The global definition is what allows the quadratic term, and later the Moreau envelope, to be applied to the whole family without boundary gluing.
\end{remark}

\section{Moreau Smoothing and the Smooth Classes}\label{app:smooth}

For a convex $m$-Lipschitz function $g:\R^n\to\R$ and $\lambda>0$ let $e_\lambda g(x)=\min_y\{g(y)+\norm{x-y}^2/(2\lambda)\}$, and let $\prox_{\lambda g}(x)$ be the unique minimizer. The following facts are standard \cite{Moreau1965,BauschkeCombettes2017}: $e_\lambda g$ is convex and differentiable, $\nabla e_\lambda g(x)=(x-\prox_{\lambda g}(x))/\lambda$, and $\nabla e_\lambda g$ is $(1/\lambda)$-Lipschitz. Moreover:
\begin{enumerate}[label=(M\arabic*),leftmargin=*]
\item $g-\lambda m^2/2\le e_\lambda g\le g$: the upper bound takes $y=x$, and the lower bound uses $g(y)\ge g(x)-m\norm{x-y}$ and $\min_{s\ge0}\{s^2/(2\lambda)-ms\}=-\lambda m^2/2$;
\item $\norm{\nabla e_\lambda g}\le m$, because $\nabla e_\lambda g(x)\in\partial g(\prox_{\lambda g}(x))$;
\item $\inf e_\lambda g=\inf g$; if $g$ attains its minimum, every minimizer of $g$ also minimizes $e_\lambda g$. Indeed, $e_\lambda g\ge\inf g$ and $e_\lambda g\le g$; at a minimizer $y^*$ this gives $e_\lambda g(y^*)=g(y^*)$.
\end{enumerate}
Isometries commute with the envelope, so $e_\lambda g_{\Theta w}(\Theta x)=e_\lambda g_w(x)$ for orthogonal $\Theta$ by Lemma~\ref{lem:belt}(i).

\begin{lemma}[Moreau hiding]\label{lem:moreau}
Let $0<a\le m$, $\rho>0$, $\theta\in(0,1/4]$, $\lambda>0$, and $e_w:=e_\lambda g_w$. There is a continuous function $H_\lambda:[0,\infty)\to\R$, independent of $w$, such that every visible $x$ satisfies
\[
|e_w(x)-H_\lambda(\norm x)|\le\min\Bigl\{a\norm x\arcsin\tfrac{|\ip wx|}{\norm x},D_1(\lambda)\Bigr\},
\]
where the first term in the minimum is read as $0$ at $x=0$, and $e_w(x)=H_\lambda(\norm x)=m\norm x-b-\lambda m^2/2$ for every visible $x$ with $\norm x\ge r_0+\lambda m$.
\end{lemma}

\begin{proof}
{Step 1: a dual formula.} Let $C:=\{(s,\mathsf p):0\le s\le1,\ \mathsf p\in smK_w\}$, a compact convex set, and $\mathsf q(s,\mathsf p):=\mathsf p-(1-s)aw$. For fixed $s$, $\max_{\mathsf p\in smK_w}\ip{\mathsf p}y=sm\sigma_w(y)$, so
\[
g_w(y)=\max_{(s,\mathsf p)\in C}\{\ip{\mathsf q(s,\mathsf p)}y-sb\}.
\]
The function $(y,(s,\mathsf p))\mapsto\ip{\mathsf q(s,\mathsf p)}y-sb+\norm{x-y}^2/(2\lambda)$ is continuous, convex in $y$, and affine in $(s,\mathsf p)$. By Sion's minimax theorem \cite{Sion1958} the minimum over $y$ and the maximum over $C$ may be interchanged, and minimizing in $y$ gives
\begin{gather*}
e_w(x)=\max_{(s,\mathsf p)\in C}\Phi_x(s,\mathsf p),\\
\Phi_x(s,\mathsf p):=\ip{\mathsf q}x-\frac\lambda2\norm{\mathsf q}^2-sb,\quad \mathsf q=\mathsf q(s,\mathsf p).
\end{gather*}

{Step 2: the gradient.} Let $(s^*,\mathsf p^*)$ maximize $\Phi_x$ over $C$ and $\mathsf q^*=\mathsf q(s^*,\mathsf p^*)$. For every $x'$, $e_w(x')\ge\Phi_{x'}(s^*,\mathsf p^*)=e_w(x)+\ip{\mathsf q^*}{x'-x}$. Thus $\mathsf q^*$ is a subgradient of the differentiable convex function $e_w$ at $x$, and $\mathsf q^*=\nabla e_w(x)$.

{Step 3: a projection.} For fixed $s^*$, $\Phi_x(s^*,\mathsf p)=-\frac\lambda2\norm{\mathsf p-(1-s^*)aw-x/\lambda}^2+\norm x^2/(2\lambda)-s^*b$, so $\mathsf p^*=\Proj_{s^*mK_w}(x/\lambda+(1-s^*)aw)$ and $\mathsf q^*=\mathsf p^*-(1-s^*)aw$.

{Step 4: tangential control.} Let $x\ne0$ be visible, and let $\kappa_x:=\min\{1/\lambda,s^*m/\norm x\}$. The projection of $x/\lambda$ onto the ball of radius $s^*m$ is $\kappa_xx$, and it belongs to $s^*mK_w$ because $x$ is visible; hence $\Proj_{s^*mK_w}(x/\lambda)=\kappa_xx$. The map $\mathrm{Id}-\Proj_{s^*mK_w}$ is firmly nonexpansive \cite{BauschkeCombettes2017}, hence nonexpansive. With $z_\lambda=x/\lambda$ and $\vartheta=(1-s^*)aw$,
\begin{align*}
\norm{\mathsf q^*-\kappa_xx}&=\bigl\|(z_\lambda+\vartheta-\Proj(z_\lambda+\vartheta))\\
&\qquad-(z_\lambda-\Proj z_\lambda)\bigr\|_2\le\norm\vartheta\le a .
\end{align*}
Consequently $|\ip{\nabla e_w(x)}\tau|=|\ip{\mathsf q^*-\kappa_xx}\tau|\le a$ for every unit vector $\tau\perp x$.

{Step 5: integration along a meridian.} Define $H_\lambda(r):=e_w(re_\perp)$ for a unit vector $e_\perp\perp w$. By the invariance under orthogonal maps, this value depends neither on $e_\perp$ nor on $w$, and $H_\lambda$ is continuous. Let $x$ be visible with $r=\norm x>0$ and $t=\ip wx$. Since $\theta<1$, $z=\norm{x-tw}>0$; put $e_\perp=(x-tw)/z$ and $\phi_0=\arcsin(t/r)$. The curve $\Gamma(\phi)=r(\cos\phi\,e_\perp+\sin\phi\,w)$ joins $re_\perp$ ($\phi=0$) to $x$ ($\phi=\phi_0$); each $\Gamma(\phi)$ is visible, has norm $r$, and $\Gamma'(\phi)\perp\Gamma(\phi)$ with $\norm{\Gamma'(\phi)}=r$. By Step 4,
\begin{align*}
|e_w(x)-H_\lambda(r)|&\le\Bigl|\int_0^{\phi_0}\ip{\nabla e_w(\Gamma(\phi))}{\Gamma'(\phi)}d\phi\Bigr|\\
&\le a r|\phi_0|\le ar\arcsin\theta .
\end{align*}

{Step 6: the exact radial tail.} Let $x$ be visible with $r\ge r_0+\lambda m$, and $y:=(1-\lambda m/r)x$. Then $y$ is visible and $\norm y\ge r_0$, so $m\sigma_w(y)-b=m\norm y-b\ge a\theta\norm y\ge-a\ip wy$: the first branch of $g_w$ is active at $y$, and every subgradient of $m\sigma_w-b$ at $y$ is a subgradient of $g_w$ at $y$. By Lemma~\ref{lem:support}, $my/\norm y\in\partial g_w(y)$. Since $(x-y)/\lambda=mx/r=my/\norm y$, the point $y$ satisfies the optimality condition of the strongly convex problem defining $e_w(x)$, and
\begin{align*}
e_w(x)&=g_w(y)+\frac{\norm{x-y}^2}{2\lambda}=m(r-\lambda m)-b+\frac{\lambda m^2}2\\
&=mr-b-\frac{\lambda m^2}2 .
\end{align*}
Equator points are visible, so the same formula holds for $H_\lambda(r)$, and the discrepancy vanishes for $r\ge r_0+\lambda m$. For $r<r_0+\lambda m$, Step 5 gives $|e_w(x)-H_\lambda(r)|\le ar\arcsin(|t|/r)$ and $ar\arcsin(|t|/r)\le a(r_0+\lambda m)\arcsin\theta=D_1(\lambda)$; for $r\ge r_0+\lambda m$ the discrepancy is $0$. At $x=0$ both sides coincide by definition.
\end{proof}

\begin{proof}[Proof of Theorem~\ref{thm:main-smooth}]
Take $m=\sqrt{L\eps}$, $\lambda=1/L$, $\rho=R$, and $a=2\eps/R$. Then $q=a/m=2\sqrt{\eps/(LR^2)}\le\frac12$ and $\upsilon:=\lambda m/\rho=\sqrt{\eps/(LR^2)}=q/2$. The function $f_w=e_w|_{B_R}$ is convex, differentiable on $\R^n$ with $L$-Lipschitz gradient, and $m$-Lipschitz by (M2), so it belongs to $\cF_{\mathrm{Sm}}^L\cap\cF_{\mathrm{Lip}}^M$ for $M\ge\sqrt{L\eps}$. Let $O_w$ equal $H_\lambda(\norm x)$ on visible points and $e_w$ on $A_w\cap B_R$, and $H(x)=H_\lambda(\norm x)$. By Lemma~\ref{lem:moreau}, $O_w\in\oo_{D_1(\lambda)}(f_w)$. By (M3) and Lemma~\ref{lem:belt}(ii), $\min_{B_R}f_w=f_w(Rw)=-2\eps$. By \eqref{eq:r0-ratio} and $1-q\theta\ge\frac78$, $D_0=2\eps\theta\frac{\theta+q}{1-q\theta}\le\frac{16}7\eps\theta(\theta+q)\le\frac37\eps$. By (M1) and Lemma~\ref{lem:belt}(iii), visible values satisfy $f_w\ge-D_0-\eps/2>-\eps$, because $D_0\le3\eps/7<\eps/2$. Their gap from the minimum $-2\eps$ is therefore strictly larger than $\eps$. Lemma~\ref{lem:hiding} yields $\Delta_{\cF}(\eps,\beta,T)\le D_1(\lambda)$ for $\cF=\cF_{\mathrm{Sm}}^L(B_R)$ and for $\cF=\cF_{\mathrm{Sm}}^L(B_R)\cap\cF_{\mathrm{Lip}}^M(B_R)$. Finally, $\arcsin\theta\le\theta/\sqrt{1-\theta^2}\le\frac{26}{25}\theta$ for $\theta\le\frac14$, $ar_0=2\eps\frac{\theta+q}{1-q\theta}\le\frac{16}{7}\eps(\theta+2\upsilon)$, and $a\lambda m=2\eps\upsilon$, so
\begin{align*}
D_1(\lambda)&\le\frac{26}{25}\eps\theta\Bigl(\frac{16}7\theta+\frac{46}{7}\upsilon\Bigr)\\
&\le\frac{12}5\eps\theta^2+7\eps\theta\upsilon .\qedhere
\end{align*}
\end{proof}

\begin{proof}[Proof of Corollary~\ref{cor:holder}]
Take $m=m_\nu$, $\lambda=\eps/m^2$, $\rho=R$, $a=2\eps/R$. Then $\lambda m^2/2=\eps/2$, $q=2\eps/(mR)\le\frac12$, and $\upsilon=\lambda m/\rho=\eps/(mR)=q/2$, so the proof of Theorem~\ref{thm:main-smooth} applies verbatim and gives $\Delta_{\cF_{\mathrm{H}}^{\nu,L_\nu}(B_R)}(\eps,\beta,T)\le\frac{12}5\eps\theta^2+7\eps^2\theta/(m_\nu R)$. The gradient of $e_w$ satisfies $\norm{\nabla e_w(x)-\nabla e_w(y)}\le\min\{\norm{x-y}/\lambda,2m\}\le(\norm{x-y}/\lambda)^\nu(2m)^{1-\nu}=2^{1-\nu}m^{1+\nu}\eps^{-\nu}\norm{x-y}^\nu=L_\nu\norm{x-y}^\nu$, so $f_w\in\cF_{\mathrm{H}}^{\nu,L_\nu}\cap\cF_{\mathrm{Lip}}^{m_\nu}$. Finally $\eps^2/m_\nu=2^{(1-\nu)/(1+\nu)}\eps^{(2+\nu)/(1+\nu)}L_\nu^{-1/(1+\nu)}$.
\end{proof}

\begin{proof}[Proof of Proposition~\ref{prop:barrier}]
Let $A_U=\E\norm U$. If $A_U=0$, then $S_\gamma g=g$, and $g(x)=M|x_1|$ has a discontinuous derivative; hence $A_U>0$. The test function $g_0=M\norm\cdot$ gives $S_\gamma g_0(0)-g_0(0)=M\gamma A_U\le\mathsf b$. If $V$ is uniform on $\Sph$ and independent of $U$, then $\E|\ip UV|=A_U\,\E|V_1|\le A_U/\sqrt n$, so some unit $v$ has $\alpha_v:=\E|\ip Uv|\le A_U/\sqrt n$. For $g_v(x)=M|\ip xv|$ the restriction $\phi(t)=S_\gamma g_v(tv)=M\E|t+\gamma Z|$, $Z=\ip Uv$, is convex with right derivative $\phi'_+(t)=M(2\Prob\{-\gamma Z\le t\}-1)$. Since $\phi'$ is $L$-Lipschitz, the distribution function of $-\gamma Z$ is $L/(2M)$-Lipschitz; hence $Z$ has a density $\mathsf f_Z$ with $\|\mathsf f_Z\|_\infty\le\gamma L/(2M)$. The existence of this bounded density excludes $\alpha_v=0$. By Markov's inequality $\Prob\{|Z|\le2\alpha_v\}\ge\frac12$, so $4\alpha_v\|\mathsf f_Z\|_\infty\ge\frac12$. Combining, $L\ge M/(4\gamma\alpha_v)\ge M\sqrt n/(4\gamma A_U)\ge M^2\sqrt n/(4\mathsf b)$.
\end{proof}

\section{Smooth Strongly Convex Class}\label{app:smsc}

\begin{lemma}[Certificate for $\cF_{\mathrm{SmSC}}^{L,\mu}$]\label{lem:smsc-general}
Let $L_0=L-\mu>0$, $\lambda=1/L_0$, $\theta$ $(n,T,\beta)$-hiding, $0<a\le m$, and $0<\rho\le R$. If
\begin{equation}\label{eq:smsc-sep}
a\rho-\frac{\mu\rho^2}{2}-D_0-\frac{\lambda m^2}{2}>\eps,
\end{equation}
then no algorithm $(\eps,\beta,D_1(\lambda),T)$-solves $\cF_{\mathrm{SmSC}}^{L,\mu}(B_R)$, and $\Delta_{\cF_{\mathrm{SmSC}}^{L,\mu}(B_R)}(\eps,\beta,T)\le D_1(\lambda)$.
\end{lemma}

\begin{proof}
Let $f_w=(e_w+Q)|_{B_R}$ with $e_w=e_\lambda g_w$. It is differentiable on $\R^n$, its gradient is $(1/\lambda+\mu)=L$-Lipschitz, and $f_w-Q=e_w$ is convex; hence $f_w\in\cF_{\mathrm{SmSC}}^{L,\mu}(B_R)$. Let $O_w$ equal $H_\lambda(\norm x)+Q(x)$ on visible points and $f_w$ on $A_w\cap B_R$. By Lemma~\ref{lem:moreau}, $O_w\in\oo_{D_1(\lambda)}(f_w)$. Since $\rho w\in B_R$, $\min_{B_R}f_w\le e_w(\rho w)+Q(\rho w)=-a\rho+\mu\rho^2/2$. Visible values satisfy $f_w\ge e_w\ge g_w-\lambda m^2/2\ge-D_0-\lambda m^2/2$ by (M1) and Lemma~\ref{lem:belt}(iii). By \eqref{eq:smsc-sep} every visible point has gap larger than $\eps$, and Lemma~\ref{lem:hiding} completes the proof.
\end{proof}

\begin{proof}[Proof of Theorem~\ref{thm:main-smsc}]
Put $\upsilon:=\sqrt{\mu/L_0}$.

(a) Take $a=3\sqrt{\mu\eps}$, $\rho=\sqrt{\eps/\mu}\le R$, $m=\sqrt{L_0\eps}$. Then $q=a/m=3\upsilon\le\frac12$, $a\rho=3\eps$, $\mu\rho^2/2=\lambda m^2/2=\eps/2$, and, by \eqref{eq:r0-ratio} and $1-q\theta\ge\frac78$, $D_0=3\eps\theta\frac{\theta+q}{1-q\theta}\le\frac{24}7\eps\theta(\theta+q)\le\frac{9}{14}\eps$, so the left side of \eqref{eq:smsc-sep} equals $2\eps-D_0>\eps$. Moreover $\lambda m/\rho=\upsilon$, $ar_0\le\frac{24}7\eps(\theta+3\upsilon)$, and $a\lambda m=3\eps\upsilon$, so
\begin{align*}
D_1(\lambda)&\le\frac{26}{25}\eps\theta\Bigl(\frac{24}7\theta+\frac{93}7\upsilon\Bigr)\\
&\le\frac{18}5\eps\theta^2+14\eps\theta\upsilon .
\end{align*}

(b) Take $a=4\sqrt{\mu\eps}$, $m=2\sqrt{L_0\eps}$, and $\rho=2\sqrt{\eps/\mu}\le R$. Then $q=a/m=2\upsilon\le2/\sqrt5<9/10$, so $a\le m$ and $1-q\theta\ge31/40$. Also $a\rho=8\eps$, $\mu\rho^2/2=2\eps$, and $\lambda m^2/2=2\eps$. Thus
\[
D_0\le\frac{320}{31}\eps\theta(\theta+q)
\le\frac{92}{31}\eps<3\eps,
\]
and the left side of \eqref{eq:smsc-sep} is $4\eps-D_0>\eps$. Moreover $ar_0\le\frac{320}{31}\eps(\theta+2\upsilon)$ and $a\lambda m=8\eps\upsilon$, hence
\begin{align*}
D_1(\lambda)&\le\frac{26}{25}\eps\theta
\left(\frac{320}{31}\theta+\frac{888}{31}\upsilon\right)\\
&\le11\eps\theta^2+30\eps\theta\upsilon .\qedhere
\end{align*}
\end{proof}

\section{The Endpoint and Its Neighborhood}\label{app:endpoint}

\begin{lemma}[Quadratic rigidity]\label{lem:quad-rigidity}
If $f$ is differentiable on a neighborhood of $B_R$, $\nabla f$ is $\mu$-Lipschitz on $B_R$, and $f-Q$ is convex on $B_R$, then $f-Q$ is affine on $B_R$.
\end{lemma}

\begin{proof}
For $x,y\in B_R$ convexity of $f-Q$ gives $f(y)\ge f(x)+\ip{\nabla f(x)}{y-x}+\frac\mu2\norm{y-x}^2$, and the $\mu$-Lipschitz gradient gives the reverse inequality by integration along the segment $[x,y]\subset B_R$. Equality holds for all $x,y$, and expanding $\norm{y-x}^2$ shows that $f(y)-\frac\mu2\norm y^2=\ip{\nabla f(x)-\mu x}{y}+f(x)-\ip{\nabla f(x)}x+\frac\mu2\norm x^2$ is affine in $y$.
\end{proof}

\begin{lemma}[Simplex reconstruction]\label{lem:simplex}
There are unit vectors $u_0,\ldots,u_n\in\R^n$ with $\sum_iu_i=0$ and $\sum_iu_iu_i^\top=\frac{n+1}{n}I$. For every $w_0,\ldots,w_n\in\R$, $\norm{\sum_iw_iu_i}^2\le\frac{n+1}n\sum_iw_i^2$.
\end{lemma}

\begin{proof}
Let $e_0,\ldots,e_n$ be the standard basis of $\R^{n+1}$ and $\mathbf 1=\sum_ie_i$. The vectors $\bar e_i=e_i-\mathbf1/(n+1)$ lie in the hyperplane $\mathbf1^\perp$, have norm $\sqrt{n/(n+1)}$, sum to $0$, and satisfy $\sum_i\bar e_i\bar e_i^\top=I-\mathbf1\mathbf1^\top/(n+1)$, the orthogonal projector onto $\mathbf1^\perp$. An isometry of $\mathbf1^\perp$ onto $\R^n$ maps $u_i:=\sqrt{(n+1)/n}\,\bar e_i$ to the required vectors. The Gram matrix of $(u_i)$ is $\frac{n+1}nI-\frac1n\mathbf1\mathbf1^\top\preceq\frac{n+1}nI$, which gives the inequality.
\end{proof}

\begin{proof}[Proof of Theorem~\ref{thm:main-endpoint}(i)--(ii)]
(i) By Lemma~\ref{lem:quad-rigidity}, $f=Q+\ip c\cdot+d'$ on $B_R$, which is $\frac\mu2\norm{x-x_c}^2+d$ with $x_c=-c/\mu$. Conversely every such function is smooth with $\mu$-Lipschitz gradient and $f-Q$ affine.

(ii) Query $y_i=O(Ru_i)$, $i=0,\ldots,n$, set $\widehat c=\frac{n}{R(n+1)}\sum_iy_iu_i$ and $\widehat x=\Proj_{B_R}(-\widehat c/\mu)$. Write $f=Q+\ip c\cdot+d'$ and $y_i=f(Ru_i)+o_i$ with $|o_i|\le\delta$. Since $Q(Ru_i)=\mu R^2/2$ is constant and $\sum_iu_i=0$, the identities of Lemma~\ref{lem:simplex} give $\widehat c-c=\frac{n}{R(n+1)}\sum_io_iu_i$ and $\norm{\widehat c-c}\le\frac{n}{R(n+1)}\sqrt{\frac{n+1}{n}(n+1)}\,\delta=\frac{\sqrt n\,\delta}{R}$. The point $\widehat x$ minimizes $\widehat f:=Q+\ip{\widehat c}\cdot$ over $B_R$. Let $x^*$ minimize $f$ over $B_R$ and $\tau_*=\norm{\widehat x-x^*}$. The optimality condition $\ip{\nabla\widehat f(\widehat x)}{x^*-\widehat x}\ge0$, the identity $\nabla f(\widehat x)=\nabla\widehat f(\widehat x)+c-\widehat c$, and $f(x^*)\ge f(\widehat x)+\ip{\nabla f(\widehat x)}{x^*-\widehat x}+\frac\mu2\tau_*^2$ give
\begin{equation}\label{eq:perturb}
f(\widehat x)-f(x^*)\le\norm{c-\widehat c}\,\tau_*-\frac\mu2\tau_*^2\le\frac{\norm{c-\widehat c}^2}{2\mu}\le\frac{n\delta^2}{2\mu R^2}.
\end{equation}
For $\delta\le R\sqrt{2\mu\eps/n}$ this is at most $\eps$ with probability one, and the algorithm, completed by arbitrary further queries when $T>n+1$, solves the class.
\end{proof}

\begin{remark}[Sharpness of the constant]\label{rem:sharp}
For odd $n$, $f=Q$, and $o_i=\delta(-1)^i$, one has $\sum_io_i=0$, $\norm{\sum_io_iu_i}^2=\frac{n+1}n(n+1)\delta^2$, and $\widehat c=\frac{n}{R(n+1)}\sum_io_iu_i$ has norm $\sqrt n\delta/R$. If $\sqrt n\delta\le\mu R^2$, then $\widehat x=-\widehat c/\mu$ and $f(\widehat x)-\min f=\norm{\widehat c}^2/(2\mu)=n\delta^2/(2\mu R^2)$, so the constant in (ii) is attained. The oracle that equals $f+o_i$ at $Ru_i$ and $f$ elsewhere is Borel.
\end{remark}

The randomized statement (iii) rests on a reflection property of adaptive linear observations of a Gaussian vector.

\begin{lemma}[Adaptive Gaussian observations]\label{lem:gauss}
Let $N\ge2$, $0\le T\le N-1$, $\zeta\sim\mathcal N(0,I_N)$, and let $\alpha_t:\R^{t-1}\to\R^N$, $t=1,\ldots,T$, be Borel maps. Define recursively $Y_t=\ip{\alpha_t(Y_1,\ldots,Y_{t-1})}{\zeta}$, $Y=(Y_1,\ldots,Y_T)$, let $\mathcal V(Y)=\spn\{\alpha_t(Y_{<t})\}$, and let $\Pi(Y)$ be the orthogonal projector onto $\mathcal V(Y)$. Let $\xi:\R^T\to\R^N$ be Borel with $\norm{\xi(Y)}=1$ and $\xi(Y)\perp\mathcal V(Y)$, and put $g=\ip{\xi(Y)}\zeta$ and $\mathcal R(\zeta)=\zeta-2g\,\xi(Y)$. Then
\begin{enumerate}[label=(\alph*),leftmargin=*]
\item the transcript of $\mathcal R(\zeta)$ equals $Y$;
\item $\mathcal R(\zeta)$ has the law $\mathcal N(0,I_N)$;
\item $g\sim\mathcal N(0,1)$.
\end{enumerate}
\end{lemma}

\begin{proof}
(a) By induction on $t$: if the first $t-1$ observations of $\mathcal R(\zeta)$ coincide with $Y_{<t}$, then the $t$-th query vector is $\alpha_t(Y_{<t})\in\mathcal V(Y)$, and $\ip{\alpha_t}{\mathcal R(\zeta)}=Y_t-2g\ip{\alpha_t}{\xi(Y)}=Y_t$.

For (b) and (c) we prove the following statement by induction on $T$, allowing the maps $\alpha_t$ to depend in a Borel way on an auxiliary random variable $\Upsilon$ independent of $\zeta$: if $\bar\zeta\sim\mathcal N(0,I_N)$ is independent of $(\Upsilon,\zeta)$, then
\begin{equation}\label{eq:gauss-rep}
(\Upsilon,Y,\zeta)\overset d=(\Upsilon,Y,\Pi(Y)\zeta+(I-\Pi(Y))\bar\zeta).
\end{equation}
Here $\Pi(Y)\zeta$ is a Borel function of $(\Upsilon,Y)$, namely the minimum-norm solution of the equations $\ip{\alpha_t}{x}=Y_t$. For $T=0$, $\Pi=0$ and \eqref{eq:gauss-rep} states $\zeta\overset d=\bar\zeta$. For the induction step, condition on $\Upsilon$ and write $u=\alpha_1/\norm{\alpha_1}$ (if $\alpha_1=0$ the first observation is identically zero and may be discarded). Then $\zeta=\ip u\zeta u+\zeta'$ with $\zeta'=(I-uu^\top)\zeta$ independent of $\ip u\zeta$, and the remaining observations are $Y_t=\ip{\alpha_t}u\ip u\zeta+\ip{(I-uu^\top)\alpha_t}{\zeta'}$, $t\ge2$. Conditionally on $\Upsilon$, the vector $u$ is fixed, $\zeta'$ is a standard Gaussian vector in $u^\perp$ independent of $Y_1=\norm{\alpha_1}\ip u\zeta$, and the quantities $\tilde Y_t=Y_t-\ip{\alpha_t}u\ip u\zeta$, $t\ge2$, form an adaptive scheme with $T-1$ observations of $\zeta'$ in which $Y_1$ is the auxiliary variable. The induction hypothesis, applied to this conditional law in $u^\perp\cong\R^{N-1}$, gives \eqref{eq:gauss-rep} for $\zeta'$ with the projector $\tilde\Pi$ onto $\tilde{\mathcal V}=\spn\{(I-uu^\top)\alpha_t\}_{t\ge2}$. Since $\mathcal V(Y)=\spn\{u\}\oplus\tilde{\mathcal V}$, $\Pi=uu^\top+\tilde\Pi$ and $I-\Pi=(I-uu^\top)-\tilde\Pi$, which yields \eqref{eq:gauss-rep} for $\zeta$.

By \eqref{eq:gauss-rep}, $g=\ip{\xi(Y)}{(I-\Pi(Y))\zeta}\overset d=\ip{\xi(Y)}{\bar\zeta}$, which is $\mathcal N(0,1)$ conditionally on $Y$; this proves (c). Similarly, since the reflection $I-2\xi\xi^\top$ commutes with $\Pi$ and fixes $\mathcal V(Y)$, $\mathcal R(\zeta)\overset d=\Pi(Y)\zeta+(I-\Pi(Y))(I-2\xi\xi^\top)\bar\zeta$, and conditionally on $Y$ the vector $(I-2\xi\xi^\top)\bar\zeta$ is standard Gaussian and independent of $Y$. Applying \eqref{eq:gauss-rep} again gives $\mathcal R(\zeta)\overset d=\zeta$, which is (b).
\end{proof}

\begin{proof}[Proof of Theorem~\ref{thm:main-endpoint}(iii)]
Fix $\eps>0$, $\beta\in(0,1/2)$, $T\le n$, and put $\eta_\beta:=\frac12(1-2\beta)>0$. For $\sigma\ge2\mu R$ let $\sigma_d=\sigma^2/\mu$, let $\zeta=(\zeta_c,\zeta_d)\sim\mathcal N(0,I_{n+1})$ with $\zeta_c\in\R^n$, and let
\[
f_\zeta(x)=Q(x)+\sigma\ip{\zeta_c}x+\sigma_d\zeta_d\in\cF_{\mathrm{SmSC}}^{\mu,\mu}(B_R).
\]
Consider first a deterministic algorithm with exact values. Its $t$-th answer is $f_\zeta(x_t)$, where $x_t$ is a Borel function of the previous answers; replacing the answers by $\tilde Y_t=f_\zeta(x_t)-Q(x_t)=\ip{(\sigma x_t,\sigma_d)}{\zeta}$ is a Borel re-encoding, so the transcript is an adaptive scheme as in Lemma~\ref{lem:gauss} with $N=n+1$ and $\alpha_t=(\sigma x_t,\sigma_d)$. For $T\ge1$ let $\mathcal D=\spn\{x_t-x_1:2\le t\le T\}$, of dimension at most $n-1$, choose a unit vector $e\perp\mathcal D$ as a Borel function of the transcript (for instance by Gram--Schmidt orthogonalization of $x_2-x_1,\ldots,x_T-x_1$ followed by the standard basis vectors), and put $k=\ip e{x_1}$, so that $\ip e{x_t}=k$ for all $t$ and $|k|\le R$; for $T=0$ let $e$ be a fixed unit vector and $k=0$. Let
\[
\varkappa:=\Bigl(1+\frac{\sigma^2k^2}{\sigma_d^2}\Bigr)^{-1/2}\in\Bigl[\frac{2}{\sqrt5},1\Bigr],\qquad \xi:=\varkappa\Bigl(e,-\frac{\sigma k}{\sigma_d}\Bigr).
\]
Then $\ip{\xi}{(\sigma x_t,\sigma_d)}=\varkappa\sigma(\ip e{x_t}-k)=0$, so $\xi\perp\mathcal V(Y)$ and Lemma~\ref{lem:gauss} applies. The reflected vector $\mathcal R(\zeta)$ has the same law and the same transcript, hence the same output $\widehat x$. Put $c=\sigma\zeta_c$, $\gamma_\sigma=\sigma\varkappa g$, and $c'=\sigma\mathcal R(\zeta)_c=c-2\gamma_\sigma e$. The success indicators $\mathcal S$ and $\mathcal S'$ of $\widehat x$ for $f_\zeta$ and $f_{\mathcal R(\zeta)}$ have equal expectations, so
\begin{equation}\label{eq:yao-half}
\Prob\{\mathcal S\}=\tfrac12\bigl(\Prob\{\mathcal S\}+\Prob\{\mathcal S'\}\bigr)\le\tfrac12\bigl(1+\Prob\{\mathcal S\cap\mathcal S'\}\bigr).
\end{equation}

We bound $\Prob\{\mathcal S\cap\mathcal S'\}$. Let $E_0=\eps+\mu R^2/2$ and $u_c=c/\norm c$ for $c\ne0$. If $\norm c\ge\mu R$, the minimizer of $f_\zeta$ over $B_R$ is $-Ru_c$, because $r\mapsto\frac\mu2r^2-\norm cr$ decreases on $[0,R]$. Since $Q(\widehat x)\ge0$, success gives $\norm c(R+\ip{u_c}{\widehat x})\le E_0$, and with $\norm{\widehat x}\le R$,
\[
\norm{u_c+\widehat x/R}^2\le2+\frac2R\ip{u_c}{\widehat x}\le\frac{2E_0}{R\norm c}.
\]
Hence on $\mathcal S\cap\mathcal S'$ with $\min\{\norm c,\norm{c'}\}\ge\mu R$, $\norm{u_c-u_{c'}}\le2\sqrt{2E_0/(R\min\{\norm c,\norm{c'}\})}$. On the other hand, with $\varsigma_e=\ip{c-\gamma_\sigma e}{e}$ we have $\ip ce=\varsigma_e+\gamma_\sigma$ and $\ip{c'}e=\varsigma_e-\gamma_\sigma$; if $|\gamma_\sigma|>|\varsigma_e|$ these have opposite signs and
\begin{align*}
\norm{u_c-u_{c'}}&\ge\frac{|\varsigma_e+\gamma_\sigma|}{\norm c}+\frac{|\varsigma_e-\gamma_\sigma|}{\norm{c'}}\\
&\ge\frac{2|\gamma_\sigma|}{\max\{\norm c,\norm{c'}\}}.
\end{align*}
A direct computation from $g=\varkappa(\ip e{\zeta_c}-\sigma k\zeta_d/\sigma_d)$ gives $\varsigma_e=\sigma^3k^2\varkappa g/\sigma_d^2+\sigma^2k\zeta_d/\sigma_d=\mu^2k^2\varkappa g/\sigma+\mu k\zeta_d$, hence $|\varsigma_e|\le\mu R(|g|/2+|\zeta_d|)$ because $|k|\le R$ and $\sigma\ge2\mu R$. For $\eta_0\in(0,1)$ and $\Lambda\ge1$ consider the event $\mathcal E=\{\norm{\zeta_c}\ge\eta_0,\ \norm{\mathcal R(\zeta)_c}\ge\eta_0,\ \norm\zeta\le\Lambda\}$. On $\mathcal E$, $\min\{\norm c,\norm{c'}\}\ge\sigma\eta_0$, $\max\{\norm c,\norm{c'}\}\le\sigma\Lambda$ (a reflection preserves $\norm\zeta$), $|g|\le\Lambda$, and $|\zeta_d|\le\Lambda$. If moreover $\sigma\eta_0\ge\mu R$, then on $\mathcal S\cap\mathcal S'\cap\mathcal E$ either $|\gamma_\sigma|\le|\varsigma_e|\le\frac32\mu R\Lambda$ or $|\gamma_\sigma|\le\sigma\Lambda\sqrt{2E_0/(R\sigma\eta_0)}$. Since $|\gamma_\sigma|\ge\sigma|g|/2$, in both cases
\[
|g|\le t_\sigma:=\frac{3\mu R\Lambda}{\sigma}+2\Lambda\sqrt{\frac{2E_0}{R\sigma\eta_0}} .
\]
By Lemma~\ref{lem:gauss}(b)--(c), $\mathcal R(\zeta)\overset d=\zeta$ and $g\sim\mathcal N(0,1)$, whose density is at most $(2\pi)^{-1/2}$, so
\[
\Prob\{\mathcal S\cap\mathcal S'\}\le2\Prob\{\norm{\zeta_c}<\eta_0\}+\Prob\{\norm\zeta>\Lambda\}+\frac{2t_\sigma}{\sqrt{2\pi}} .
\]
Choose $\eta_0$ and $\Lambda$ such that the first two terms are at most $\eta_\beta/2$, and then $\sigma$ so large that $\sigma\eta_0\ge\mu R$ and $2t_\sigma/\sqrt{2\pi}\le\eta_\beta/2$. These choices do not depend on the algorithm, and \eqref{eq:yao-half} gives $\Prob\{\mathcal S\}\le\frac12(1+\eta_\beta)<1-\beta$ for every deterministic algorithm with exact values.

A randomized algorithm is a Borel mixture over $\omega$ of deterministic algorithms, and its success indicator is jointly Borel in $(\omega,\zeta)$ because the maps $\chi_t$ and $(\zeta,x)\mapsto f_\zeta(x)$ are Borel. By Fubini's theorem, $\E_\zeta\Prob_A\{\text{success on }f_\zeta\}=\E_\omega\Prob_\zeta\{\text{success}\}<1-\beta$, so some $f_\zeta$ is solved with probability less than $1-\beta$ even by the exact oracle $O=f_\zeta\in\oo_0(f_\zeta)$.
\end{proof}

\begin{proof}[Proof of Theorem~\ref{thm:main-endpoint}(iv)]
Let $a=\sqrt{3\mu\eps}$ and $f_w(x)=Q(x)-a\ip wx$. Then $f_w\in\cF_{\mathrm{SmSC}}^{L,\mu}(B_R)$ for every $L\ge\mu$, its minimizer $(a/\mu)w$ lies in $B_R$ because $\eps\le\mu R^2/3$, and $\min_{B_R}f_w=-a^2/(2\mu)=-\frac32\eps$. For visible $x$, $f_w(x)\ge\frac\mu2\norm x^2-a\theta\norm x\ge-a^2\theta^2/(2\mu)=-\frac32\eps\theta^2\ge-\frac3{32}\eps$, so every visible point has gap strictly larger than $\eps$, and every $\eps$-minimizer lies in $A_w$. The oracle equal to $Q$ on visible points and $f_w$ on $A_w\cap B_R$ has error at most $a\theta R$. Lemma~\ref{lem:hiding} with $H=Q$ gives the claim, and $\theta=\theta_*$ gives the last display of the theorem.
\end{proof}

\begin{lemma}[Structure of $\cF_{\mathrm{SmSC}}^{L,\mu}$]\label{lem:smsc-structure}
If $f\in\cF_{\mathrm{SmSC}}^{L,\mu}(B_R)$, then $L\ge\mu$ and $\tilde f:=f-Q$ is convex with $(L-\mu)$-Lipschitz gradient on $B_R$.
\end{lemma}

\begin{proof}
On the open ball $\mathring B_R=\operatorname{int}B_R$ the gradient $\nabla f$ is Lipschitz, hence differentiable almost everywhere by Rademacher's theorem \cite{EvansGariepy2015}, with a Jacobian $J(x)$ that coincides almost everywhere with the distributional Hessian of $f$ and is therefore symmetric. Since $\nabla f$ is $L$-Lipschitz, $\norm{J(x)}\le L$. Since $\tilde f$ is convex and differentiable, $x\mapsto\nabla\tilde f(x)=\nabla f(x)-\mu x$ is monotone, so $v^\top(J(x)-\mu I)v\ge0$ for all $v$. The symmetric matrix $J(x)-\mu I$ therefore has spectrum in $[0,L-\mu]$, which forces $L\ge\mu$ and gives $\norm{J(x)-\mu I}\le L-\mu$ almost everywhere. For $x,y\in\mathring B_R$ and almost every small translate $(x',y')$ of $(x,y)$, the segment $[x',y']\subset\mathring B_R$ meets the exceptional set in a set of length zero by Fubini's theorem; integrating the Lipschitz map $\nabla\tilde f$ along it gives $\norm{\nabla\tilde f(y')-\nabla\tilde f(x')}\le(L-\mu)\norm{y'-x'}$. Continuity of $\nabla\tilde f$ extends this to $x,y$ and then to $B_R$.
\end{proof}

\begin{proof}[Proof of Proposition~\ref{prop:near-end}]
Let $c=\nabla\tilde f(0)$ and $\varrho(x)=\tilde f(x)-\tilde f(0)-\ip cx$. By Lemma~\ref{lem:smsc-structure}, $0\le\varrho(x)\le\frac{L_0}2\norm x^2$ and $\norm{\nabla\varrho(x)}\le L_0\norm x$ on $B_R$, and $f=Q+\ip c\cdot+\tilde f(0)+\varrho$. Write $O(ru_i)=f(ru_i)+o_i$ with $|o_i|\le\delta$. Since $Q(ru_i)$ and $\tilde f(0)$ do not depend on $i$ and $\sum_iu_i=0$,
\[
\widehat c-c=\frac{n}{r(n+1)}\sum_i\Bigl(\varrho(ru_i)-\frac{L_0r^2}4+o_i\Bigr)u_i ,
\]
and every coefficient has modulus at most $\frac{L_0r^2}4+\delta$. Lemma~\ref{lem:simplex} gives $\norm{\widehat c-c}\le\sqrt n\bigl(\frac{L_0r}4+\frac\delta r\bigr)$. The argument leading to \eqref{eq:perturb}, now with $\nabla f(\widehat x)=\nabla\widehat f(\widehat x)+(c-\widehat c)+\nabla\varrho(\widehat x)$ and $f$ $\mu$-strongly convex, gives $f(\widehat x)-\min_{B_R}f\le(\norm{\widehat c-c}+L_0R)^2/(2\mu)$, which is \eqref{eq:near-end}. The choice $r=R$ gives the first consequence. If $0<\delta\le L_0R^2/4$, the choice $r=2\sqrt{\delta/L_0}\le R$ gives $\frac{L_0r}4+\frac\delta r=\sqrt{L_0\delta}$, and the right side of \eqref{eq:near-end} is at most $\eps$ whenever $\delta\le(\sqrt{2\mu\eps}-L_0R)^2/(nL_0)$; this gives the second consequence. The exact oracle is covered by monotonicity from any positive noise level satisfying these bounds.
\end{proof}

\section{Achievability, Budgets, and Effective Dimension}\label{app:ach}

Let $h\in(0,R)$, $R'=R-h$, and let $f$ be convex and $M$-Lipschitz on $B_R$, with the convex $M$-Lipschitz extension $\bar f(x)=\inf_{y\in B_R}\{f(y)+M\norm{x-y}\}$ to $\R^n$. Let $U$ be uniform on the unit ball, $e$ uniform on $\Sph$, and $f_h(x)=\E\bar f(x+hU)$.

\begin{lemma}[Ball smoothing]\label{lem:ballsmooth}
The function $f_h$ is convex and differentiable on $\R^n$ with $\nabla f_h(x)=\frac n{2h}\E[(\bar f(x+he)-\bar f(x-he))e]$. On $B_{R'}$, $f\le f_h\le f+Mh$, and if $f-Q$ is convex on $B_R$, then $f_h-Q$ is convex on $B_{R'}$.
\end{lemma}

\begin{proof}
The formula $\nabla f_h(x)=\frac nh\E[\bar f(x+he)e]$ for continuous $\bar f$ is \cite[Lemma~1]{FlaxmanKalaiMcMahan2005}; replacing $e$ by $-e$ gives the symmetric form. Convexity is preserved by averaging. For $x\in B_{R'}$ the points $x+hU$ lie in $B_R$, where $\bar f=f$; Jensen's inequality and $\E U=0$ give $f_h(x)\ge f(x)$, and the Lipschitz bound gives $f_h\le f+Mh$. Finally $\E Q(x+hU)=Q(x)+\frac\mu2h^2\E\norm U^2$, so $f_h-Q$ is an average of the convex functions $f(\cdot+hU)-Q(\cdot+hU)$ on $B_{R'}$ plus a constant.
\end{proof}

Let $O\in\oo_\delta(f)$ and, for $x\in B_{R'}$, $\widehat G(x,e)=\frac n{2h}(O(x+he)-O(x-he))e$. By Lemma~\ref{lem:ballsmooth},
\begin{equation}\label{eq:bias}
\begin{aligned}
\norm{\E_e\widehat G(x,e)-\nabla f_h(x)}&\le\frac{n}{2h}\sup_{\norm u=1}\E\bigl[2\delta|\ip eu|\bigr]\\
&\le\frac{\sqrt n\,\delta}h=:\mathsf B,
\end{aligned}
\end{equation}
because $\E|\ip eu|\le(\E\ip eu^2)^{1/2}=n^{-1/2}$, and $\norm{\widehat G}^2\le n^2(M+\delta/h)^2=:S_G$ deterministically. Consider $x_{k+1}=\Proj_{B_{R'}}(x_k-\eta_k\widehat G(x_k,e_k))$ with independent $e_k$ and $x_1=0$, and let $x^*$ minimize $f_h$ over $B_{R'}$.

\begin{lemma}[Projected two-point method]\label{lem:sgd}
Let $S\ge S_G$. (a) With $\eta_k=R'/\sqrt{KS}$ and $\bar x=\frac1K\sum_{k=1}^Kx_k$, $\E f_h(\bar x)-f_h(x^*)\le R'\sqrt{S/K}+2R'\mathsf B$.
(b) If $f-Q$ is convex on $B_R$, then with $\eta_k=4/(\mu(k+1))$ and $\bar x=\sum_{k=1}^Kkx_k/\sum_{k=1}^Kk$, $\E f_h(\bar x)-f_h(x^*)\le\frac{4S}{\mu(K+1)}+\frac{\mathsf B^2}{\mu}$.
\end{lemma}

\begin{proof}
Let $d_k=\E\norm{x_k-x^*}^2$ and $b_k=\E[\widehat G(x_k,e_k)\mid x_k]-\nabla f_h(x_k)$. Nonexpansiveness of the projection gives $\norm{x_{k+1}-x^*}^2\le\norm{x_k-x^*}^2-2\eta_k\ip{\widehat G}{x_k-x^*}+\eta_k^2S$. (a) By convexity, $\ip{\nabla f_h(x_k)}{x_k-x^*}\ge f_h(x_k)-f_h(x^*)$, and $|\ip{b_k}{x_k-x^*}|\le2R'\mathsf B$. Summing over $k$ and using $d_1\le R'^2$ and Jensen's inequality gives $\E f_h(\bar x)-f_h(x^*)\le\frac{R'^2}{2\eta K}+\frac{\eta S}2+2R'\mathsf B$, and the choice of $\eta$ gives (a). (b) By Lemma~\ref{lem:ballsmooth}, $f_h$ is $\mu$-strongly convex on $B_{R'}$, so $\ip{\nabla f_h(x_k)}{x_k-x^*}\ge f_h(x_k)-f_h(x^*)+\frac\mu2\norm{x_k-x^*}^2$, and $\ip{b_k}{x_k-x^*}\ge-\frac\mu4\norm{x_k-x^*}^2-\frac{\mathsf B^2}\mu$. Hence
\begin{align*}
\E f_h(x_k)-f_h(x^*)&\le\frac{\mu(k-1)}8d_k-\frac{\mu(k+1)}8d_{k+1}\\
&\quad+\frac{2S}{\mu(k+1)}+\frac{\mathsf B^2}\mu .
\end{align*}
Multiplying by $k$ and summing, the first two terms telescope to a nonpositive quantity, and division by $\sum_kk=K(K+1)/2$ together with Jensen's inequality gives (b).
\end{proof}

\begin{proof}[Proof of Propositions~\ref{prop:ach-lip} and \ref{prop:ach-sc}]
Take $h=7\eps/(64M)$ in the Lipschitz case and $h=\eps/(8M)$ in the strongly convex case; then $h\le R/8$ because $\eps\le MR$ (respectively $\eps\le\mu R^2\le MR$, since $M\ge\mu R$ by Proposition~\ref{prop:lipsc-end}). Let $x_f$ minimize $f$ over $B_R$. Then $\min_{B_{R'}}f_h\le f_h(\Proj_{B_{R'}}x_f)\le f(x_f)+2Mh$, and $f\le f_h$ on $B_{R'}$ by Lemma~\ref{lem:ballsmooth}, so every output $\bar x\in B_{R'}$ satisfies $\E f(\bar x)-\min_{B_R}f\le\E f_h(\bar x)-f_h(x^*)+2Mh$.

Lipschitz case: $2Mh=\frac7{32}\eps$. By \eqref{eq:bias}, the bound $\delta\le\eps^2/(84\sqrt nMR)$ gives\[2R'\mathsf B\le\frac{2R\sqrt n\delta}h=\frac{128\sqrt nMR\delta}{7\eps}\le\frac{32}{147}\eps\le\frac7{32}\eps\]and $\delta/h\le\frac{16\eps}{147\sqrt nR}\le\frac{16}{147\sqrt2}M$, because $\eps\le MR$ and $n\ge2$. Since $(1+16/(147\sqrt2))^2\le29/25$, we may take $S=(29/25)n^2M^2$. The stated $K$ then gives $R'\sqrt{S/K}\le\sqrt{29/30000}\,\eps\le\eps/32$. Lemma~\ref{lem:sgd}(a) with this $S$ yields $\E f(\bar x)-\min f\le(\frac1{32}+\frac7{32}+\frac7{32})\eps=\frac{15}{32}\eps$.

Strongly convex case: $2Mh=\frac14\eps$. By \eqref{eq:bias}, $\delta\le\sqrt\mu\eps^{3/2}/(19\sqrt nM)$ gives $\mathsf B^2/\mu\le\frac{64}{361}\eps$ and $\delta/h\le8\sqrt{\mu\eps}/(19\sqrt n)\le\frac{8}{19\sqrt2}M<3M/10$, because $\sqrt{\mu\eps}\le\mu R\le M$. Hence we may take $S=(169/100)n^2M^2$, and the stated $K$ gives $4S/(\mu(K+1))\le\eps/32$. Lemma~\ref{lem:sgd}(b) yields
\[
\E f(\bar x)-\min f\le\left(\frac1{32}+\frac{64}{361}+\frac14\right)\eps<\frac{23}{50}\eps.
\]

Moreover $\delta\le\eps/(84\sqrt2)<\eps/32$ in the Lipschitz case and $\delta\le\eps\sqrt{\mu\eps}/(19\sqrt nM)\le\eps/(19\sqrt2)<\eps/25$ in the strongly convex case. By Markov's inequality a single run of $2K$ queries returns a point with gap at most $g_1:=\frac{15}{16}\eps$ in the Lipschitz case and $g_1:=23\eps/25$ in the strongly convex case with probability at least $\frac12$. Run $J$ independent copies, query $O$ at the $J$ outputs, and return the output with the smallest oracle value. With probability at least $1-2^{-J}\ge1-\beta$ some run succeeds, and the returned point $\widehat x$ then satisfies $f(\widehat x)\le O(\widehat x)+\delta\le O(\bar x^{(j)})+\delta\le f(\bar x^{(j)})+2\delta$ for the successful run $j$, hence $f(\widehat x)-\min f\le g_1+2\delta$, which is at most $\frac{15}{16}\eps+\frac{1}{16}\eps=\eps$ and $(23/25+2/25)\eps=\eps$, respectively. The total number of queries is $J(2K+1)$. All maps involved are Borel, so this is a randomized algorithm in the sense of Section~\ref{sec:model}.
\end{proof}

\begin{proof}[Proof of Proposition~\ref{prop:ach-smsc}]
Let $f\in\cF_{\mathrm{SmSC}}^{L,\mu,\circ}(B_R)$ and let $x_f$ be its minimizer over $B_R$. Since $\norm{x_f}<R$, $\nabla f(x_f)=0$, so $\norm{\nabla f(x)}\le L\norm{x-x_f}\le2LR$ on $B_R$, and $f$ is $2LR$-Lipschitz; the constructions above apply with $M=2LR$. Take $h=\sqrt{5\eps/(16L)}$, so that $h<\sqrt{\eps/L}\le R$ and $x_f\in B_{R'}$. For $x\in B_{R'}$ and $\norm u\le1$ the segment $[x,x+hu]$ lies in $B_R$, so $f(x+hu)\le f(x)+h\ip{\nabla f(x)}u+\frac L2h^2$; averaging with $\E U=0$ gives $f_h\le f+\frac L2h^2=f+\frac5{32}\eps$ on $B_{R'}$. Together with $f\le f_h$ on $B_{R'}$ (Lemma~\ref{lem:ballsmooth}) and $\min_{B_{R'}}f_h\le f_h(x_f)$, every output $\bar x\in B_{R'}$ satisfies $\E f(\bar x)-\min_{B_R}f\le\E f_h(\bar x)-f_h(x^*)+\frac5{32}\eps$.

For $\delta\le\frac\eps5\sqrt{\mu/(nL)}$, \eqref{eq:bias} gives $\mathsf B^2/\mu=n\delta^2/(\mu h^2)\le\frac{16}{125}\eps$, and $\delta/h\le\frac4{5\sqrt5}\sqrt{\mu\eps/n}<\frac12LR$ because $\sqrt{\mu\eps}\le\mu R\le LR$; hence $S_G\le(25/4)n^2L^2R^2=:S$, and $K\ge800n^2L^2R^2/(\mu\eps)$ gives $\frac{4S}{\mu(K+1)}\le\frac\eps{32}$. Lemma~\ref{lem:sgd}(b) with this $S$, applicable since $f-Q$ is convex on $B_R$, yields $\E f(\bar x)-\min f\le(\frac5{32}+\frac{16}{125}+\frac1{32})\eps<\eps/3$. By Markov's inequality a single run of $2K$ queries returns a point with gap at most $2\eps/3$ with probability at least $\frac12$, and the selection among $J$ runs described in the proof of Propositions~\ref{prop:ach-lip} and \ref{prop:ach-sc} returns, with probability at least $1-\beta$, a point with gap at most $2\eps/3+2\delta\le2\eps/3+\eps/3=\eps$, since $\delta\le\eps/(5\sqrt2)<\eps/6$.
\end{proof}

\begin{proof}[Proof of Corollary~\ref{cor:two-sided}]
Part (a) combines Theorem~\ref{thm:main-lip} and Proposition~\ref{prop:ach-lip}. For (b), the hard functions of Theorem~\ref{thm:main-smooth} lie in $\cF_{\mathrm{Sm}}^L\cap\cF_{\mathrm{Lip}}^M$ for $M\ge\sqrt{L\eps}$, which gives the upper bound; the class is contained in $\cF_{\mathrm{Lip}}^M$, and $\eps\le\kappa\sqrt{L\eps}R$ holds since $\eps\le LR^2/16$, so Proposition~\ref{prop:ach-lip} gives the lower bound $\eps^2/(84\sqrt nMR)$ with $M=\kappa\sqrt{L\eps}$. For (c), the assumed inequality $M-\mu R\ge4\sqrt{\mu\eps}$ and $M\ge2\mu R$ give $M-\mu R\ge M/2$, so Theorem~\ref{thm:main-lipsc} applies with $(M-\mu R)^{-1}\le2M^{-1}$, and Proposition~\ref{prop:ach-sc} gives the lower bound. For (d), the hard functions of Corollary~\ref{cor:holder} lie in $\cF_{\mathrm{H}}^{\nu,L_\nu}\cap\cF_{\mathrm{Lip}}^{m_\nu}\subset\cF$, the class $\cF$ is contained in $\cF_{\mathrm{Lip}}^M$, and $\eps\le m_\nu R\le MR$, so Proposition~\ref{prop:ach-lip} gives the lower bound with $M=\kappa m_\nu$. For (e), let $f_w=(e_w+Q)|_{B_R}$ be the hard functions of the proof of Theorem~\ref{thm:main-smsc}(b), with $a=4\sqrt{\mu\eps}$, $m=2\sqrt{L_0\eps}$, $\rho=2\sqrt{\eps/\mu}$, $\lambda=1/L_0$, and let $x_w$ minimize $f_w$ over $B_R$. By (M3) and Lemma~\ref{lem:belt}(ii), $f_w(x)\ge-a\rho+\frac\mu2\norm x^2$ for every $x$, and $f_w(x_w)\le f_w(\rho w)\le-a\rho+\frac\mu2\rho^2$. Hence $\norm{x_w}\le\rho=2\sqrt{\eps/\mu}$. Since $L\ge6\mu$ and $\eps\le\mu R^2/6$, $\norm{x_w}+\sqrt{\eps/L}\le(2+1/\sqrt6)\sqrt{\eps/\mu}<\sqrt6\sqrt{\eps/\mu}\le R$, so $f_w\in\cF_{\mathrm{SmSC}}^{L,\mu,\circ}(B_R)$. Lemma~\ref{lem:hiding} applied to this subclass, with the oracles of the proof of Lemma~\ref{lem:smsc-general}, gives the upper bound, and Proposition~\ref{prop:ach-smsc} gives the lower bound.
\end{proof}

\begin{proof}[Proof of Corollary~\ref{cor:budget}]
For $\theta=\theta(\delta)$, $\frac43\eps\theta^2+\frac53\frac{\eps^2\theta}{MR}\le\frac43\cdot\frac{3\delta}8+\frac53\cdot\frac{3\delta}{10}=\delta$, so an algorithm that $(\eps,\beta,\delta,T)$-solves $\cF_{\mathrm{Lip}}^M(B_R)$ also solves it at the level $D_0\le\delta$ of the proof of Theorem~\ref{thm:main-lip}. By Lemma~\ref{lem:hiding}, $\theta(\delta)$ is therefore not hiding: $2(T+1)(1-\theta^2)^{(n-1)/2}\ge1-\beta$, and $(1-\theta^2)^{(n-1)/2}\le e^{-(n-1)\theta^2/2}$ gives the claim.
\end{proof}

\begin{proof}[Proof of Corollary~\ref{cor:effdim}]
Let $w$ be uniform on the unit sphere of $V$, let $g_w$ be the belt function of the same parameters, and replace every hard function $f_w$ of Theorems~\ref{thm:main-lip}--\ref{thm:main-smsc} and Corollary~\ref{cor:holder} by $f_w^V:=(f_w-Q_\cF)\circ P+Q_\cF$. For the smooth classes, $f_w^V-Q_\cF=e_\lambda g_w\circ P$, and Lemma~\ref{lem:moreau} applies at $Px$ because $\ip w{Px}=\ip wx$; composition with $P$ preserves convexity, Lipschitz constants, and the Lipschitz or H\"older constants of gradients, and $f_w^V-Q_\cF$ depends only on $Px$. The hidden sets are $A^V_w=\{x:Px\ne0,\ |\ip wx|>\theta\norm{Px}\}$, the common oracle is the corresponding radial function of $Px$ plus $Q_\cF$, and Lemmas~\ref{lem:belt} and \ref{lem:moreau}, applied at $Px$, give the same discrepancies at every point $x$ through $Px$. The comparison points $\rho w$ lie in $V\cap B_R$ and retain the objective values used in the original separation arguments. They need not minimize the strongly convex smoothed instances. The additional term $Q_\cF(x)-Q_\cF(Px)$ is nonnegative, so the lower bounds on visible objective values are preserved. For fixed $x$ with $Px\ne0$, Lemma~\ref{lem:cap} in $V$ gives $\Prob\{x\in A^V_w\}\le2(1-\theta^2)^{(d-1)/2}$, and points with $Px=0$ are visible. Lemma~\ref{lem:hiding} then applies with $n$ replaced by $d$.
\end{proof}

\EOD

\begin{thebibliography}{99}

\bibitem{NemirovskiYudin1983}
A.~S. Nemirovsky and D.~B. Yudin, \emph{Problem Complexity and Method Efficiency in Optimization}. Transl. by E.~R. Dawson. Chichester, U.K.: Wiley, 1983. Available: \url{https://www2.isye.gatech.edu/~nemirovs/Nemirovskii_Yudin_1983.pdf}.

\bibitem{ConnScheinbergVicente2009}
A.~R. Conn, K.~Scheinberg, and L.~N. Vicente, \emph{Introduction to Derivative-Free Optimization}. Philadelphia, PA, USA: SIAM, 2009.
doi: \url{https://doi.org/10.1137/1.9780898718768}.

\bibitem{Spall2003}
J.~C. Spall, \emph{Introduction to Stochastic Search and Optimization: Estimation, Simulation, and Control}. Hoboken, NJ, USA: Wiley, 2003.
doi: \url{https://doi.org/10.1002/0471722138}.

\bibitem{LarsonMenickellyWild2019}
J.~Larson, M.~Menickelly, and S.~M. Wild, ``Derivative-free optimization methods,'' \emph{Acta Numerica}, vol.~28, pp.~287--404, 2019.
doi: \url{https://doi.org/10.1017/S0962492919000060}.

\bibitem{FlaxmanKalaiMcMahan2005}
A.~D. Flaxman, A.~T. Kalai, and H.~B. McMahan, ``Online convex optimization in the bandit setting: gradient descent without a gradient,'' in \emph{Proc. 16th Annual ACM-SIAM Symposium on Discrete Algorithms (SODA)}, 2005, pp.~385--394. Available: \url{https://research.google/pubs/online-convex-optimization-in-the-bandit-setting-gradient-descent-without-a-gradient/}.

\bibitem{BubeckLeeEldan2017}
S.~Bubeck, Y.~T. Lee, and R.~Eldan, ``Kernel-based methods for bandit convex optimization,'' in \emph{Proc. 49th Annual ACM Symposium on Theory of Computing (STOC)}, 2017, pp.~72--85.
doi: \url{https://doi.org/10.1145/3055399.3055403}.

\bibitem{DevolderGlineurNesterov2014}
O.~Devolder, F.~Glineur, and Y.~Nesterov, ``First-order methods of smooth convex optimization with inexact oracle,'' \emph{Mathematical Programming}, vol.~146, nos.~1--2, pp.~37--75, 2014.
doi: \url{https://doi.org/10.1007/s10107-013-0677-5}.

\bibitem{dAspremont2008}
A.~d'Aspremont, ``Smooth optimization with approximate gradient,'' \emph{SIAM Journal on Optimization}, vol.~19, no.~3, pp.~1171--1183, 2008.
doi: \url{https://doi.org/10.1137/060676386}.

\bibitem{BerahasCaoChoromanskiScheinberg2022}
A.~S. Berahas, L.~Cao, K.~Choromanski, and K.~Scheinberg, ``A theoretical and empirical comparison of gradient approximations in derivative-free optimization,'' \emph{Foundations of Computational Mathematics}, vol.~22, no.~2, pp.~507--560, 2022.
doi: \url{https://doi.org/10.1007/s10208-021-09513-z}.

\bibitem{Protasov1996}
V.~Yu. Protasov, ``Algorithms for approximate calculation of the minimum of a convex function from its values,'' \emph{Mathematical Notes}, vol.~59, no.~1, pp.~69--74, 1996.
doi: \url{https://doi.org/10.1007/BF02312467}.

\bibitem{DvinskikhTominin2022}
D.~Dvinskikh, V.~Tominin, I.~Tominin, and A.~Gasnikov,
``Noisy zeroth-order optimization for non-smooth saddle point problems,''
in \emph{Mathematical Optimization Theory and Operations Research
(MOTOR 2022)}, Lecture Notes in Computer Science, vol.~13367.
Cham, Switzerland: Springer, 2022, pp.~18--33.
doi: \url{https://doi.org/10.1007/978-3-031-09607-5_2}.

\bibitem{AlashqarGasnikovDvinskikhLobanov2023}
B.~A. Alashqar, A.~V. Gasnikov, D.~M. Dvinskikh, and A.~V. Lobanov,
``Gradient-free federated learning methods with $l_1$ and
$l_2$-randomization for non-smooth convex stochastic optimization
problems,''
\emph{Computational Mathematics and Mathematical Physics},
vol.~63, no.~9, pp.~1600--1653, 2023.
doi: \url{https://doi.org/10.1134/S0965542523090026}.

\bibitem{BelloniLiangNarayananRakhlin2015}
A.~Belloni, T.~Liang, H.~Narayanan, and A.~Rakhlin, ``Escaping the local minima via simulated annealing: Optimization of approximately convex functions,'' in \emph{Proc. Conference on Learning Theory}, PMLR, vol.~40, 2015, pp.~240--265. Available: \url{https://proceedings.mlr.press/v40/Belloni15.html}.

\bibitem{SingerVondrak2015}
Y.~Singer and J.~Vondr\'ak, ``Information-theoretic lower bounds for convex optimization with erroneous oracles,'' in \emph{Advances in Neural Information Processing Systems}, vol.~28, 2015, pp.~3204--3212. Available: \url{https://papers.nips.cc/paper_files/paper/2015/hash/393c55aea738548df743a186d15f3bef-Abstract.html}.

\bibitem{LiRisteski2016}
Y.~Li and A.~Risteski, ``Algorithms and matching lower bounds for approximately-convex optimization,'' in \emph{Advances in Neural Information Processing Systems}, vol.~29, 2016, pp.~4745--4753. Available: \url{https://papers.nips.cc/paper_files/paper/2016/hash/186fb23a33995d91ce3c2212189178c8-Abstract.html}.

\bibitem{NesterovSpokoiny2017}
Y.~Nesterov and V.~Spokoiny, ``Random gradient-free minimization of convex functions,'' \emph{Foundations of Computational Mathematics}, vol.~17, no.~2, pp.~527--566, 2017.
doi: \url{https://doi.org/10.1007/s10208-015-9296-2}.

\bibitem{Gasnikov2024}
A.~Gasnikov, D.~Dvinskikh, P.~Dvurechensky, E.~Gorbunov, A.~Beznosikov, and A.~Lobanov, ``Randomized gradient-free methods in convex optimization,'' in \emph{Encyclopedia of Optimization}, P.~M. Pardalos and O.~A. Prokopyev, Eds. Cham, Switzerland: Springer, 2024, pp.~1--15. doi: \url{https://doi.org/10.1007/978-3-030-54621-2_859-1}.

\bibitem{GorbunovDvurechenskyGasnikov2022}
E.~Gorbunov, P.~Dvurechensky, and A.~Gasnikov, ``An accelerated method for derivative-free smooth stochastic convex optimization,'' \emph{SIAM Journal on Optimization}, vol.~32, no.~2, pp.~1210--1238, 2022. doi: \url{https://doi.org/10.1137/19M1259225}.

\bibitem{DvurechenskyGorbunovGasnikov2021}
P.~Dvurechensky, E.~Gorbunov, and A.~Gasnikov, ``An accelerated directional derivative method for smooth stochastic convex optimization,'' \emph{European Journal of Operational Research}, vol.~290, no.~2, pp.~601--621, 2021.
doi: \url{https://doi.org/10.1016/j.ejor.2020.08.027}.

\bibitem{BogolubskyEtAl2016}
L.~Bogolubsky, P.~Dvurechensky, A.~Gasnikov, G.~Gusev, Y.~Nesterov, A.~Raigorodskii, A.~Tikhonov, and M.~Zhukovskii, ``Learning supervised PageRank with gradient-based and gradient-free optimization methods,'' in \emph{Advances in Neural Information Processing Systems}, vol.~29, 2016. Available: \url{https://papers.nips.cc/paper_files/paper/2016/hash/1f34004ebcb05f9acda6016d5cc52d5e-Abstract.html}.

\bibitem{Kornilov2023GradientFree}
N.~Kornilov, A.~Gasnikov, P.~Dvurechensky, and D.~Dvinskikh, ``Gradient-free methods for non-smooth convex stochastic optimization with heavy-tailed noise on convex compact,'' \emph{Computational Management Science}, vol.~20, Art.~no.~37, 2023. doi: \url{https://doi.org/10.1007/s10287-023-00470-2}.

\bibitem{Kornilov2023Heavy}
N.~Kornilov, O.~Shamir, A.~Lobanov, D.~Dvinskikh, A.~Gasnikov, I.~Shibaev, E.~Gorbunov, and S.~Horv\'ath, ``Accelerated zeroth-order method for non-smooth stochastic convex optimization problem with infinite variance,'' in \emph{Advances in Neural Information Processing Systems}, vol.~36, 2023, pp.~64083--64102.
doi: \url{https://doi.org/10.52202/075280-2800}.

\bibitem{DuchiJordanWainwrightWibisono2015}
J.~C. Duchi, M.~I. Jordan, M.~J. Wainwright, and A.~Wibisono, ``Optimal rates for zero-order convex optimization: The power of two function evaluations,'' \emph{IEEE Transactions on Information Theory}, vol.~61, no.~5, pp.~2788--2806, 2015.
doi: \url{https://doi.org/10.1109/TIT.2015.2409256}.

\bibitem{Shamir2017}
O.~Shamir, ``An optimal algorithm for bandit and zero-order convex optimization with two-point feedback,'' \emph{Journal of Machine Learning Research}, vol.~18, no.~52, pp.~1--11, 2017. Available: \url{https://jmlr.org/papers/v18/16-632.html}.

\bibitem{BayandinaGasnikovLagunovskaya2018}
A.~S. Bayandina, A.~V. Gasnikov, and A.~A. Lagunovskaya,
``Gradient-free two-point methods for solving stochastic nonsmooth
convex optimization problems with small non-random noises,''
\emph{Automation and Remote Control},
vol.~79, no.~8, pp.~1399--1408, 2018.
doi: \url{https://doi.org/10.1134/S0005117918080039}.

\bibitem{LobanovBashirovGasnikov2024}
A.~Lobanov, N.~Bashirov, and A.~Gasnikov,
``The `black-box' optimization problem: Zero-order accelerated
stochastic method via kernel approximation,''
\emph{Journal of Optimization Theory and Applications},
vol.~203, no.~3, pp.~2451--2486, 2024.
doi: \url{https://doi.org/10.1007/s10957-024-02556-6}.

\bibitem{GasnikovLobanovStonyakin2024}
A.~V. Gasnikov, A.~V. Lobanov, and F.~S. Stonyakin,
``Highly smooth zeroth-order methods for solving optimization
problems under the PL condition,''
\emph{Computational Mathematics and Mathematical Physics},
vol.~64, no.~4, pp.~739--770, 2024.
doi: \url{https://doi.org/10.1134/S0965542524700118}.

\bibitem{Lobanov2023Stochastic}
A.~Lobanov,
``Stochastic adversarial noise in the `black box' optimization problem,''
in \emph{Optimization and Applications (OPTIMA 2023)},
Lecture Notes in Computer Science, vol.~14395.
Cham, Switzerland: Springer, 2023, pp.~60--71.
doi: \url{https://doi.org/10.1007/978-3-031-47859-8_5}.

\bibitem{LobanovGasnikov2024LipschitzNoise}
A.~Lobanov and A.~Gasnikov, ``Improved maximum noise level estimation in black-box optimization problems,'' \emph{Zap. Nauchn. Sem. POMI}, vol.~540, pp.~132--147, 2024. Available: \url{https://www.mathnet.ru/eng/znsl7547}.

\bibitem{BachPerchet2016}
F.~Bach and V.~Perchet, ``Highly-smooth zero-th order online optimization,'' in \emph{Proc. Conference on Learning Theory}, PMLR, vol.~49, 2016, pp.~257--283. Available: \url{https://proceedings.mlr.press/v49/bach16.html}.

\bibitem{AkhavanPontilTsybakov2020}
A.~Akhavan, M.~Pontil, and A.~B. Tsybakov, ``Exploiting higher order smoothness in derivative-free optimization and continuous bandits,'' in \emph{Advances in Neural Information Processing Systems}, vol.~33, 2020, pp.~9017--9027. Available: \url{https://papers.nips.cc/paper_files/paper/2020/hash/6646b06b90bd13dabc11ddba01270d23-Abstract.html}.

\bibitem{AkhavanPontilTsybakov2021}
A.~Akhavan, M.~Pontil, and A.~B. Tsybakov, ``Distributed zero-order optimization under adversarial noise,'' in \emph{Advances in Neural Information Processing Systems}, vol.~34, 2021, pp.~10209--10220. Available: \url{https://papers.nips.cc/paper_files/paper/2021/hash/5487e79fa0ccd0b79e5d4a4c8ced005d-Abstract.html}.

\bibitem{AkhavanChzhenPontilTsybakov2024}
A.~Akhavan, E.~Chzhen, M.~Pontil, and A.~B. Tsybakov, ``Gradient-free optimization of highly smooth functions: Improved analysis and a new algorithm,'' \emph{Journal of Machine Learning Research}, vol.~25, no.~370, pp.~1--50, 2024. Available: \url{https://www.jmlr.org/papers/v25/23-0733.html}.

\bibitem{BychkovEtAl2024}
G.~K. Bychkov, D.~M. Dvinskikh, A.~V. Antsiferova, A.~V. Gasnikov, and A.~V. Lobanov, ``Accelerated zero-order SGD under high-order smoothness and overparameterized regime,'' \emph{Russian Journal of Nonlinear Dynamics}, vol.~20, no.~5, pp.~759--788, 2024. doi: \url{https://doi.org/10.20537/nd241209}.

\bibitem{JamiesonNowakRecht2012}
K.~G. Jamieson, R.~Nowak, and B.~Recht, ``Query complexity of derivative-free optimization,'' in \emph{Advances in Neural Information Processing Systems}, vol.~25, 2012, pp.~2672--2680. Available: \url{https://papers.nips.cc/paper_files/paper/2012/hash/e6d8545daa42d5ced125a4bf747b3688-Abstract.html}.

\bibitem{Shamir2013}
O.~Shamir, ``On the complexity of bandit and derivative-free stochastic convex optimization,'' in \emph{Proc. Conference on Learning Theory}, PMLR, vol.~30, 2013, pp.~3--24. Available: \url{https://proceedings.mlr.press/v30/Shamir13.html}.

\bibitem{Kerger2026}
P.~Kerger, ``Closing the oracle-complexity gap in derivative-free convex optimization: A near-quadratic lower bound from exact function values,'' arXiv:2607.13335, 2026.
doi: \url{https://doi.org/10.48550/arXiv.2607.13335}.

\bibitem{ZhangZhangQiLin2026}
H.~Zhang, C.~Zhang, Z.~Qi, and Z.~Lin, ``Near-optimal lower bounds for randomized algorithms in exact value zeroth-order convex optimization,'' arXiv:2607.16558, 2026.
doi: \url{https://doi.org/10.48550/arXiv.2607.16558}.

\bibitem{Nesterov2015Universal}
Y.~Nesterov, ``Universal gradient methods for convex optimization problems,'' \emph{Mathematical Programming}, vol.~152, nos.~1--2, pp.~381--404, 2015.
doi: \url{https://doi.org/10.1007/s10107-014-0790-0}.

\bibitem{IouditskiNesterov2014}
A.~Iouditski and Y.~Nesterov, ``Primal-dual subgradient methods for minimizing uniformly convex functions,'' arXiv:1401.1792, 2014.
doi: \url{https://doi.org/10.48550/arXiv.1401.1792}.

\bibitem{DuchiBartlettWainwright2012}
J.~C. Duchi, P.~L. Bartlett, and M.~J. Wainwright, ``Randomized smoothing for stochastic optimization,'' \emph{SIAM Journal on Optimization}, vol.~22, no.~2, pp.~674--701, 2012.
doi: \url{https://doi.org/10.1137/110831659}.

\bibitem{Nesterov2008Cubic}
Y.~Nesterov, ``Accelerating the cubic regularization of Newton's method on convex problems,'' \emph{Mathematical Programming}, vol.~112, no.~1, pp.~159--181, 2008.
doi: \url{https://doi.org/10.1007/s10107-006-0089-x}.

\bibitem{BauschkeCombettes2017}
H.~H. Bauschke and P.~L. Combettes, \emph{Convex Analysis and Monotone Operator Theory in Hilbert Spaces}, 2nd~ed. Cham, Switzerland: Springer, 2017.
doi: \url{https://doi.org/10.1007/978-3-319-48311-5}.

\bibitem{Moreau1965}
J.-J. Moreau, ``Proximit\'e et dualit\'e dans un espace hilbertien,'' \emph{Bulletin de la Soci\'et\'e Math\'ematique de France}, vol.~93, pp.~273--299, 1965.
doi: \url{https://doi.org/10.24033/bsmf.1625}.

\bibitem{Sion1958}
M.~Sion, ``On general minimax theorems,'' \emph{Pacific Journal of Mathematics}, vol.~8, no.~1, pp.~171--176, 1958.
doi: \url{https://doi.org/10.2140/pjm.1958.8.171}.

\bibitem{EvansGariepy2015}
L.~C. Evans and R.~F. Gariepy, \emph{Measure Theory and Fine Properties of Functions}, revised ed. Boca Raton, FL, USA: CRC Press, 2015.
doi: \url{https://doi.org/10.1201/b18333}.

\end{thebibliography}
\end{document}